\catcode`\@=11

\magnification=1200
\hsize=13.5cm     \vsize=21cm
\hoffset=0mm    \voffset=-5mm
\pretolerance=500  \tolerance=1000  \brokenpenalty=5000

\baselineskip=15pt
\parskip=5pt plus 2pt minus 1pt

\catcode`\;\active
\def;{\relax\ifhmode\ifdim\lastskip>\z@
\unskip\fi\kern.2em\fi\string;}

\catcode`\:=\active
\def:{\relax\ifhmode\ifdim\lastskip>\z@
\unskip\fi\penalty\@M\ \fi\string:}

\catcode`\!\active
\def!{\relax\ifhmode\ifdim\lastskip>\z@
\unskip\fi\kern.2em\fi\string!}

\catcode`\?\active
\def?{\relax\ifhmode\ifdim\lastskip>\z@
\unskip\fi\kern.2em\fi\string?}

\def\^#1{\if#1i{\accent"5E\i}\else{\accent"5E #1}\fi}
\def\"#1{\if#1i{\accent"7F\i}\else{\accent"7F #1}\fi}

\frenchspacing

\newcount\notenumber  \notenumber=1
\def\note#1{\footnote{($^{\the\notenumber}$)}{#1}\global%
\advance\notenumber by 1}

\font\pcap=cmcsc10
\font\bfit=cmbxti10
\font\twelvebf=cmbx12

\def\cf{{\it cf.\/}\ }  \def\ie{{\it i.e.\/}\ }
\def\eg{{\it e.g.\/}\ }

\def\og{\leavevmode\raise
.3ex\hbox{$\scriptscriptstyle\langle\!\langle\,$}}
\def\fg{\leavevmode\raise
.3ex\hbox{$\scriptscriptstyle\,\rangle\!\rangle$}}

\newif\ifpagetitre   \pagetitretrue
\newtoks\hautpagetitre   	 \hautpagetitre={\hfill}
\newtoks\baspagetitre \baspagetitre={\hfill}

\newtoks\auteur   	\auteur={\hfill}
\newtoks\titrecourant   	\titrecourant={\hfill}

\newtoks\hautpagegauche   \newtoks\hautpagedroite
\hautpagegauche={\hfill}  \hautpagedroite={\hfill}

\newtoks\baspagegauche    \newtoks\baspagedroite
\baspagegauche={\hfill\tenrm --\ \folio\ --\hfill}
\baspagedroite={\hfill\tenrm --\ \folio\ --\hfill}

\headline={\ifpagetitre\the\hautpagetitre
\else\ifodd\pageno
\the\hautpagedroite\else\the\hautpagegauche\fi\fi}

\footline={\ifpagetitre\the\baspagetitre\global\pagetitrefalse
\else\ifodd\pageno
\the\baspagedroite\else\the\baspagegauche\fi\fi}

\newtoks\proof    \proof={D\'emonstration}

\def\pf{\par \vskip-\medskipamount\vskip-2pt\noindent{\it \the\proof.}\/ }

\def\cqfd{\vskip-1mm\nobreak{\line{\hfill\vbox{\hrule height 4pt
width 4pt}\quad}}\medskip\penalty -300}

\def\ra{\rightarrow}  \def\Ra{\Rightarrow}
\def\sd{\hbox{$\times\!$\vbox{\hrule height 5pt width .5pt}}\,}
\def\dst{\displaystyle}

\def\C{{\bf C}}
\def\R{{\bf R}}

\def\N{{\bf N}}
\def\L{{\cal L}}
\def\F{{\cal F}}
\def\mt{\mapsto}

\def\id{{\rm id}}

\catcode`\@=12

\newcount\thno  \thno=1
\newcount\bibno  \bibno=1
\newcount\secno  \secno=0
\newcount\exerno  \exerno=1

\def\sn{\the\secno .}

\def\bn{\item{[\the\bibno ]}\global\advance\bibno by 1}

\def\sec#1{\global\advance\secno by 1
\global\thno=1\global \exerno=1
\vskip 18pt plus 15pt\penalty -1000\centerline{\twelvebf\sn\
{#1}}\nobreak}

\newcount\svbibno\svbibno=1

\def\savebib#1{\expandafter\xdef\csname XA#1\endcsname%
{\global\svbibno=\number\bibno\relax}%
\global\expandafter\def\csname #1\endcsname%
{\csname XA#1\endcsname\the\svbibno}}

\def\labbib#1{\savebib{#1}\global\advance\bibno by 1\relax}

\def\exer{\medskip\noindent\the \exerno .\ \global\advance \exerno by 1}

\newtoks\theorem         \theorem={Th\'eor\`eme}
\newtoks\proposition    \proposition={Proposition}
\newtoks\definition      \definition={D\'efinition}
\newtoks\lemma            \lemma={Lemme}
\newtoks\corollary      \corollary={Corollaire}
\newtoks\remark         \remark={Remarque}
\newtoks\remarks        \remarks={Remarques}
\newtoks\example        \example={Exemple}
\newtoks\examples      \examples={Exemples}
\newtoks\theorema      \theorema={le th\'eor\`eme}
\newtoks\propositiona  \propositiona={la proposition}
\newtoks\definitiona   \definitiona={la d\'efinition}
\newtoks\lemmaa       \lemmaa={le lemme}
\newtoks\corollarya   \corollarya={le corollaire}
\newtoks\exercise  \exercise={Exercice}
\newtoks\exercises \exercises={Exercices}

\def\exerc{\medskip\noindent{\it \the\exercise \the\exerno .}\/ ---\
\global\advance\exerno by 1}

\def\exercs{\medskip\noindent{\the\exercises.}\/ ---\
\the\exerno .\ \global\advance\exerno by 1}
\def\th#1{\vskip 12pt plus 9 pt\penalty-300\noindent\sn\the\thno .\ {\pcap
\the\theorem}. ---\ {\sl {#1}}\global\advance\thno by 1
\penalty-300\medskip}

\newcount\svsecno  \svsecno=0
\newcount\svthno\svthno=1

\def\thl#1#2{{\expandafter\xdef\csname XA#1ns\endcsname%
{\global\svthno=\number\thno\relax%
\global\svsecno=\number\secno\relax}%
\global\expandafter\def\csname #1ns\endcsname%
{\csname XA#1ns\endcsname\the\svthno}}
{\expandafter\xdef\csname XA#1s\endcsname%
{\global\svthno=\number\thno\relax%
\global\svsecno=\number\secno\relax}%
\global\expandafter\def\csname #1s\endcsname%
{\csname XA#1s\endcsname\the\svsecno.\the\svthno}}
{\expandafter\xdef\csname XA#1\endcsname%
{\global\svthno=\number\thno\relax%
\global\svsecno=\number\secno\relax}%
\global\expandafter\def\csname #1\endcsname%
{\csname XA#1\endcsname\the\theorema\ \the\svsecno.\the\svthno}}
\medskip\penalty-200\noindent\sn\the\thno
.\ {\pcap \the\theorem}. ---\ {\sl {#2}}\global\advance\thno by 1
\penalty-200\medskip}

\def\df#1{\medskip\penalty-200\noindent\sn\the\thno .\ {\pcap
\the\definition}. ---\ {\sl {#1}} \global\advance\thno
by 1\medskip\penalty-200}

\def\dfl#1#2{{\expandafter\xdef\csname XA#1ns\endcsname%
{\global\svthno=\number\thno\relax%
\global\svsecno=\number\secno\relax}%
\global\expandafter\def\csname #1ns\endcsname%
{\csname XA#1ns\endcsname\the\svthno}}
{\expandafter\xdef\csname XA#1s\endcsname%
{\global\svthno=\number\thno\relax%
\global\svsecno=\number\secno\relax}%
\global\expandafter\def\csname #1s\endcsname%
{\csname XA#1s\endcsname\the\svsecno.\the\svthno}}
{\expandafter\xdef\csname XA#1\endcsname%
{\global\svthno=\number\thno\relax%
\global\svsecno=\number\secno\relax}%
\global\expandafter\def\csname #1\endcsname%
{\csname XA#1\endcsname\the\definitiona\ \the\svsecno.\the\svthno}}
\medskip\penalty-200
\noindent\sn\the\thno .\ {\pcap \the\definition}. ---\ {\sl {#2}}
\global\advance\thno by 1\medskip\penalty-200}

\def\cor#1{\medskip\penalty-200\noindent\sn\the\thno .\ {\pcap
\the\corollary}. ---\ {\sl {#1}}\global\advance\thno by 1
\medskip\penalty-200}

\def\corl#1#2{{\expandafter\xdef\csname XA#1ns\endcsname%
{\global\svthno=\number\thno\relax%
\global\svsecno=\number\secno\relax}%
\global\expandafter\def\csname #1ns\endcsname%
{\csname XA#1ns\endcsname\the\svthno}}
{\expandafter\xdef\csname XA#1s\endcsname%
{\global\svthno=\number\thno\relax%
\global\svsecno=\number\secno\relax}%
\global\expandafter\def\csname #1s\endcsname%
{\csname XA#1s\endcsname\the\svsecno.\the\svthno}}
{\expandafter\xdef\csname XA#1\endcsname%
{\global\svthno=\number\thno\relax%
\global\svsecno=\number\secno\relax}%
\global\expandafter\def\csname #1\endcsname%
{\csname XA#1\endcsname\the\corollarya\ \the\svsecno.\the\svthno}}
\medskip\penalty-200
\noindent\sn\the\thno .\ {\pcap \the\corollary}. ---\ {\sl {#2}}
\global\advance\thno by 1\medskip\penalty-200}

\def\prop#1{\vskip 12pt plus 9 pt\penalty-200\noindent\sn\the\thno .\
{\pcap Proposition}. ---\ {\sl {#1}}\global\advance\thno by 1
\medskip\penalty-200}

\def\propl#1#2{{\expandafter\xdef\csname XA#1ns\endcsname%
{\global\svthno=\number\thno\relax%
\global\svsecno=\number\secno\relax}%
\global\expandafter\def\csname #1ns\endcsname%
{\csname XA#1ns\endcsname\the\svthno}}
{\expandafter\xdef\csname XA#1s\endcsname%
{\global\svthno=\number\thno\relax%
\global\svsecno=\number\secno\relax}%
\global\expandafter\def\csname #1s\endcsname%
{\csname XA#1s\endcsname\the\svsecno.\the\svthno}}
{\expandafter\xdef\csname XA#1\endcsname%
{\global\svthno=\number\thno\relax%
\global\svsecno=\number\secno\relax}%
\global\expandafter\def\csname #1\endcsname%
{\csname XA#1\endcsname\the\propositiona\ \the\svsecno.\the\svthno}}
\medskip\penalty-200
\noindent\sn\the\thno .\ {\pcap \the\proposition}. ---\ {\sl {#2}}
\global\advance\thno by 1\medskip\penalty-200}

\def\lem#1{\vskip 12pt plus 9 pt\penalty-200\noindent\sn\the\thno .\
{\pcap
\the\lemma}. ---\ {\sl {#1}}\global\advance\thno by 1
\medskip\penalty-200}

\def\leml#1#2{{\expandafter\xdef\csname XA#1ns\endcsname%
{\global\svthno=\number\thno\relax%
\global\svsecno=\number\secno\relax}%
\global\expandafter\def\csname #1ns\endcsname%
{\csname XA#1ns\endcsname\the\svthno}}
{\expandafter\xdef\csname XA#1s\endcsname%
{\global\svthno=\number\thno\relax%
\global\svsecno=\number\secno\relax}%
\global\expandafter\def\csname #1s\endcsname%
{\csname XA#1s\endcsname\the\svsecno.\the\svthno}}
{\expandafter\xdef\csname XA#1\endcsname%
{\global\svthno=\number\thno\relax%
\global\svsecno=\number\secno\relax}%
\global\expandafter\def\csname #1\endcsname%
{\csname XA#1\endcsname\the\lemmaa\ \the\svsecno.\the\svthno}}
\medskip\penalty-200
\noindent\sn\the\thno .\ {\pcap \the\lemma}. ---\ {\sl {#2}}
\global\advance\thno by 1\medskip\penalty-200}

\def\rem{\medskip\penalty-200\noindent{\it\sn\the\thno .
\the\remark.}\/ ---\ \global\advance\thno by 1}

\def\rems{\medskip\penalty-200\noindent{\it\sn\the\thno .
\the\remarks.}\/ ---\ \global\advance\thno by 1}

\def\ex{\medskip\penalty-200\noindent{\it\sn\the\thno .
\the \example.}\/ ---\ \global\advance\thno by 1}

\def\exs{\medskip\penalty-200\noindent{\it\sn\the\thno .
\the \examples.}\/ ---\ \global\advance\thno by 1}

\def\subsec{\vskip 12pt plus 9 pt\penalty-200\noindent \sn\the\thno .
---\ \global\advance\thno by 1}
\def\H{{\cal H}}
\def\ssec#1{\bigskip\noindent{\bfit #1}}

\def\ref{{}}

\labbib {BSb}
\labbib {Da}
\labbib{Dab}
\labbib{Enock}
\labbib{Enockb}
\labbib {EN}
\labbib{ES}
\labbib{GHJ}
\labbib{HS}
\labbib{IK}
\labbib {ILP}
\labbib {Kac}
\labbib{KP}
\labbib {Longo}
\labbib {Maj}
\labbib {Ma}
\labbib {Mon}
\labbib {NZ}
\labbib{Po}
\labbib {Sz}
\labbib{Takeu}
\labbib{Wa}

\font\Bigtm=cmr10 at 20pt

\hfill\vfill\vfill
\centerline{\Bigtm Unitaires multiplicatifs en dimension finie}
\bigskip
\centerline{\Bigtm et leurs sous-objets.}
\vfill

\centerline{\bf Saad BAAJ, Etienne BLANCHARD et Georges
SKANDALIS.}

\vfill
\vfill\vfill\vfill

\parindent 0pt\baselineskip13pt

SB, D\'epartement de Math\'ematiques, Universit\'e Blaise
Pascal, F-63177 Aubi\`ere.\hfill\break
e-mail: baaj@ucfma.univ-bpclermont.fr.

EB, Institut de Math\'ematiques de Luminy, Case
907, F-13288 Marseille Cedex 9.\hfill\break
e-mail: blanch@iml.univ-mrs.fr

GS, Institut de Math\'ematiques de Jussieu,
Universit\'e Denis Diderot (Paris 7), {\pcap Case Postale} {\bf 7012}, 2,
Place Jussieu, F-75251 PARIS C\'edex 05.\hfill\break
e-mail: skandal@math.jussieu.fr.

\vfill\eject
\hfill
\bigskip\bigskip\bigskip\bigskip
\centerline{\twelvebf Unitaires multiplicatifs en dimension finie et leurs
sous-objets.}

\bigskip\bigskip\bigskip\centerline{Saad {\pcap Baaj}, Etienne {\pcap
Blanchard} et Georges {\pcap Skandalis}}

\vfill\vfill\vfill

\noindent {\bf R\'esum\'e: }On appelle {\it pr\'e-sous-groupe}\/ d'un unitaire
multiplicatif $\;V\;$ agissant sur un espace hilbertien de dimension finie
$\;\H\;$ une droite vectorielle $\;L\;$ de $\;\H\;$ telle que $\;V(L\otimes
L)=L\otimes L\,$. Nous montrons que les pr\'e-sous-groupes sont en nombre
fini, donnons un \'equivalent du th\'eor\`eme de Lagrange et
g\'en\'eralisons \`a ce cadre la construction du `bi-produit crois\'e'. De plus,
nous \'etablissons des bijections entre pr\'e-sous-groupes et
sous-alg\`ebres co\"id\'eales de l'alg\`ebre de Hopf associ\'ee \`a $\;V\,$,
et donc avec les facteurs interm\'ediaires des inclusions de facteurs
associ\'ees (\cf [\ILP]). Enfin, nous montrons que les pr\'e-sous-groupes
classifient les sous-objets de $\;(\H,V)\,$.

\bigskip\bigskip 
\noindent {\bf Abstract: }A {\it pre-subgroup}\/ of a
multiplicative unitary
$\;V\;$ on a finite dimensionnal Hilbert space $\;\H\;$ is a vector line
$\;L\;$ in $\;\H\;$ such that $\;V(L\otimes L)=L\otimes L\,$. We show
that there are finitely many  pre-subgroups, give a Lagrange theorem and
generalize the construction of a `bi-crossed product'. Moreover,
we establish bijections between pre-subgroups and coideal subalgebras
of the Hopf algebra associated with $\;V\,$, and therefore with the
intermediate subfactors of the associated (depth two) inclusions (\cf
[\ILP]). Finally, we show that the pre-subgroups classify the subobjects
of $\;(\H,V)\,$.

\vfill\vfill\vfill\eject
\centerline{\twelvebf Unitaires multiplicatifs en dimension finie et leurs
sous-objets.}

\bigskip\centerline{Saad {\pcap Baaj}, Etienne {\pcap Blanchard} et Georges
{\pcap Skandalis}}

\bigskip\bigskip Dans cet article, nous \'etudions l'\'equivalent des
sous-groupes pour les alg\`ebres de Kac de dimension finie. Nous montrons
que les \og sous-groupes\fg\ en notre sens correspondent aux
sous-alg\`ebres co\"id\'eales des alg\`ebres de Hopf (\cf \eg [\ref\Mon]) et
donc aux facteurs interm\'ediaires des inclusions de profondeur 2 (\cf
[\ref\ILP]).

\bigskip Le point de vue pris ici est celui de [\ref\BSb]: nous partons d'un
unitaire multiplicatif $\;V\in{\cal L}({\cal H}\otimes {\cal H})\;$
o\`u $\;{\cal H}\;$ est un espace hilbertien de dimension finie. Une question
tr\`es naturelle est alors: {\sl Quels sont les sous-objets de $\;({\cal
H},V)\,$?} Autrement dit, quels sont les sous-espaces $\;H\;$ de $\;{\cal
H}\;$ tels que $\;V(H\otimes H)=H\otimes H\,$?

Cette question nous am\`ene \`a \'etudier d'abord le cas des sous-espaces
comme ci-dessus de dimension 1, c'est \`a dire les vecteurs
unitaires $\;f\in{\cal H}\;$ tels que $\;V(f\otimes f)=f\otimes f\,$. Dans le
cas de l'unitaire multiplicatif d'un groupe fini $\;G\,$, ($\,{\cal
H}=\ell^2(G)\,$, $\;V(\xi)(s,t)=\xi (st,t)\,,$ pour $\;\xi \in {\cal H}\otimes
{\cal H}=\ell^2(G\times G)\,,\;s,t\in G\,$) un tel $\;f\;$ est (proportionnel
\`a) la fonction caract\'eristique d'un sous-groupe de $\;G\,$.

Dans le cas g\'en\'eral, nous appelons un tel $\;f\;$ un {\it
pr\'e-sous-groupe}\/ de $\;V\,$. A un tel pr\'e-sous-groupe, correspondent
deux sous-espaces de $\;{\cal H}\,$, not\'es $\;H^f\;$ et $\;H_f\,$, qui dans
le cas d'un sous-groupe $\;\Gamma\;$ d'un groupe $\;G\;$ consistent
respectivement en le sous-espace des fonctions nulles hors de
$\;\Gamma\;$ (l'espace $\;H^f\,$) et invariantes par translation \`a droite
par $\;\Gamma\;$ (l'espace $\;H_f\,$).

Nous obtenons ais\'ement un \'equivalent dans notre cadre du th\'eor\`eme de
Lagrange:  pour tout pr\'e-sous-groupe $\;f\,,$ on a $\;\dim H^f\dim
H_f=\dim{\cal H}\,$.

Il y a naturellement une relation d'ordre $\;\prec\;$ pour les
pr\'e-sous-groupes: l'inclusion pour les $\;H^f\;$ correspondants. Cette
relation a un plus grand \'el\'ement, \og le\fg\ vecteur fixe $\;e\;$ et un
plus-petit \'el\'ement, \og le\fg\ vecteur co-fixe $\;\hat e\,$. De plus,
l'espace des pr\'e-sous-groupes est un treillis. Cela nous permet de d\'efinir
l'analogue dans notre cadre du sous-groupe engendr\'e par une partie.

Si $\;f\;$ et $\;g\;$ sont deux pr\'e-sous-groupes,  $\;|\langle
f,g\rangle|^{-2} \;$ est un entier. Il en r\'esulte qu'il y a un nombre fini de
pr\'e-sous-groupes et on peut m\^eme donner une estimation (assez
grossi\`ere) de leur nombre.

\medskip Nous disons qu'un pr\'e-sous-groupe $\;f\;$ est un {\it
sous-groupe}\/ (resp. un {\it co-sous-groupe}) si $\;V(H^f\otimes
H^f)=H^f\otimes H^f\;$ (resp. $\;V(H_f\otimes H_f)=H_f\otimes H_f\,$);
nous montrons que cela a lieu si et seulement si le projecteur sur $\;H^f\;$
(resp. $\;H_f\,$) est dans le centre de l'alg\`ebre $\;S\;$ (resp. $\;\widehat
S\,$) associ\'ee \`a $\;V\;$\note{Pour l'unitaire multiplicatif associ\'e \`a
un groupe fini $\;G\,$, l'alg\`ebre $\;S=C(G)\;$ associ\'ee est commutative,
donc tout pr\'e-sous-groupe est un sous-groupe.}. Si $\;f\;$ est \`a la fois un
sous-groupe et un co-sous-groupe, on dit que c'est un {\it sous-groupe
normal}. Les notions de sous-groupe, co-sous-groupe et sous-groupe normal
co\"incident clairement avec celles introduites par Kac dans [\ref\Kac]. Les
sous-groupes et les co-sous-groupes forment un sous-treillis de l'ensemble
des pr\'e-sous-groupes.

Enfin, les pr\'e-sous-groupes permettent de classifier les {\it
sous-quotients}\/  de $\;V\;$ \ie les sous-espaces $\;H\;$ de $\;{\cal H}\;$
tels que $\;V(H\otimes H)=H\otimes H\,$: pour un tel $\;H\,$, il existe un
unique couple $\;(f,\hat f)\;$ de pr\'e-sous-groupes tels que $\;\hat f\prec
f\;$ et $\;H=H^f\cap H_{\hat f}\,$. Inversement, nous donnons plusieurs
conditions n\'ecessaires et suffisantes pour qu'une telle intersection soit un
sous-quotient. En particulier, nous montrons que pour tout
pr\'e-sous-groupe $\;\hat f\;$ il existe un plus grand pr\'e-sous-groupe
$\;f\;$ (le {\it normalisateur}\/ de $\;\hat f\,$) tel que $\;\hat f\prec f\;$
et $\;H^f\cap H_{\hat f}\;$ soit un sous-quotient.

Si $\;f\;$ est un pr\'e-sous-groupe, l'ensemble $\;D_f=\{\,x\in S\,,\;
\kappa(x)e\in H_f\,\}\;$ (o\`u $\;\kappa \;$ d\'esigne l'antipode de $\;S\,$)
est une sous-alg\`ebre co\"id\'eale (\`a droite) \ie une sous-alg\`ebre
involutive de $\;S\;$ telle que $\;\delta (D_f)\subset D_f\otimes S\;$ (o\`u
$\;\delta \;$ d\'esigne le coproduit de $\;S\,$). R\'eciproquement, toute
sous-alg\`ebre co\"id\'eale (\`a droite) $\;D\;$ de $\;S\;$ est de la forme
$\;D_f\,$. En ce sens, notre analogue du th\'eor\`eme de Lagrange est une
r\'e\'ecriture d'un r\'esultat plus g\'en\'eral de [\ref\Ma]. Notons que dans
[\ref\Ma], ce th\'eor\`eme de Lagrange r\'esulte d'un r\'esultat plus pr\'ecis:
$\;S\;$ est un $\,D$-module libre, dont nous donnons, dans notre cadre, une
d\'emonstration nouvelle. Nous discutons aussi de nouveaux r\'esultats de
libert\'e analogues \`a ceux de [\ref\Ma]: plus pr\'ecis\'ement, nous
d\'efinissons aussi une sous-alg\`ebre co\"id\'eale (\`a gauche)
$\;\widehat G_f=\{\;y\in \widehat S\,,\; \widehat \kappa (y)\hat e\in
H^f\,\}\;$ de $\;\widehat S\,$. Nous montrons que l'espace vectoriel
$\;B_{f,f}\;$ engendr\'e par $\;\{\,xy\,,\;x\in D_f\,,\;y\in \widehat
G_f\,\}\;$ est une sous-alg\`ebre involutive de $\;{\cal L}({\cal H})\;$ et que
$\;{\cal H}\;$ est un module libre (de rang 1) sur $\;B_{f,f}\,$. Nous
g\'en\'eralisons cela \`a une alg\`ebre $\;B_{f,g}\;$ associ\'ee de fa\c con
analogue \`a deux pr\'e-sous-groupes $\;f\;$ et $\;g\,$.

\bigskip Le coproduit de l'alg\`ebre de Kac $\;S\;$ (et donc l'unitaire
multiplicatif $\;V\,$) d\'epend uniquement des alg\`ebres $\;S\;$ et
$\;\widehat S\;$ et des vecteurs $\;e\;$ et $\;\hat e\,$. On peut en fait
caract\'eriser les quadruplets $\;(S,\widehat S,e,\hat e)\;$ parmi les
quadruplets $\;(A,B,f,\hat f)\;$ o\`u $\;A\;$ et $\;B\;$ sont des
sous-$C^*$-alg\`ebres de $\;\L(\H)\;$ et $\;f,\hat f\;$ sont des vecteurs (de
norme 1) de $\;\H\,$.

A l'aide de cette caract\'erisation nous donnons une g\'en\'eralisation \`a
notre cadre du biproduit crois\'e (\cf [\Kac, \Takeu, \Maj]): \`a deux
pr\'e-sous-groupes $\;f\;$ et $\;g\;$  les plus \'eloign\'es possible \ie
tels que $\;|\langle f,g\rangle |^{-2}=\dim {\cal H}\,$, nous associons un
nouvel unitaire multiplicatif $\;W\,$; l'alg\`ebre $\;\widehat S\;$ associ\'ee
\`a $\;W\;$ est l'alg\`ebre $\;B_{g,g}\;$ associ\'ee \`a $\;g\,$;
l'alg\`ebre $\;S\;$ est une alg\`ebre $\;A_{f,f}\;$ construite de fa\c con
analogue. Cette construction s'interpr\`ete simplement en termes
d'inclusions de facteurs.

Dans le cas de deux sous-groupes $\;G_1\,,\;G_2\;$ d'un groupe fini $\;G\;$
tels que $\;G_1G_2=G\,,\;G_1\cap G_2=\{e\}\,$, l'unitaire multiplicatif
$\;W\;$ est associ\'e au \og biproduit crois\'e\fg\ de [\ref\Takeu,\Maj] (\cf
aussi [\Kac] et [\ref\BSb] \S8).

\medskip De plus, par [\ref\ILP], les sous-alg\`ebres co\"id\'eales
correspondent aux facteurs interm\'ediaires d'une inclusion irr\'eductible
associ\'ee \`a $\;S\,$. Plusieurs de nos r\'esultats s'interpr\`etent en
termes de facteurs interm\'ediaires. En particulier, notre r\'esultat de
finitude est une r\'e\'ecriture (plus pr\'ecise) d'un r\'esultat plus
g\'en\'eral
de [\ref\Wa].

\bigskip L'organisation de ce travail est la suivante:

{\parindent 0pt Dans le premier paragraphe, nous faisons quelques rappels et
compl\'ements pour les unitaires multiplicatifs sur un espace hilbertien de
dimension finie.

Dans le deuxi\`eme paragraphe, nous donnons une caract\'erisation des
alg\`ebres de Kac (et leur alg\`ebre duale) parmi les couples de
sous-alg\`ebres de $\;\L(\H)\,$.

Nous introduisons ensuite les pr\'e-sous-groupes, pour lesquels nous
donnons plusieurs r\'esultats d'int\'egralit\'e, ainsi que la relation d'ordre
et ses propri\'et\'es (paragraphe 3). Nous faisons le lien avec les
sous-alg\`ebres co\"id\'eales au quatri\`eme paragraphe, ainsi que la
construction de plusieurs autres sous-alg\`ebres associ\'ees \`a des
pr\'e-sous-groupes.

Au cinqui\`eme paragraphe, nous discutons les sous-groupes,
co-sous-groupes et sous-quotients.

Au sixi\`eme paragraphe, en nous aidant du paragraphe 2, nous donnons la
g\'en\'eralisation du biproduit crois\'e mentionn\'ee ci-dessus.

Nous interpr\'etons enfin (paragraphe 7) nos r\'esultats en termes
d'inclusions de facteurs d'indice fini et de profondeur 2.

}

\bigskip Nous remercions M. Izumi qui a attir\'e notre attention sur le lien
entre sous-alg\`ebres co\"id\'eales et sous-facteurs de [\ILP] et D. Bisch
pour plusieurs conseils bibliographiques.

\sec{Pr\'eliminaires.}

Dans ce paragraphe, nous fixons quelques notations et particularisons
\`a la dimension finie certaines constructions de [\ES\ref] et
[\BSb\ref]. Nous donnons en outre quelques autres r\'esultats plus
sp\'ecifiques \`a la dimension finie (\cf aussi [\KP]).

\noindent {\bf Notations}. --- Dans cet article, $\;{\cal H}\;$
d\'esigne un espace hilbertien de dimension finie $\;n\ne 0\;$ et
$\;V\in {\cal L}({\cal H}\otimes {\cal H})\;$ un unitaire multiplicatif
sur $\;{\cal H}\,$, \ie tel que $\;V_{12}V_{13}V_{23}=V_{23}V_{12}\
(\in {\cal L}({\cal H} \otimes {\cal H}\otimes {\cal H})\,$).

Rappelons (\cf [\ref\BSb]) que pour $\;\omega \in {\cal L}({\cal
H})^*\;$ on note $\;L(\omega )=(\omega\otimes \id )(V)\in {\cal
L}({\cal H})\;$ et $\;\rho(\omega )=(\id\otimes \omega )(V)\,$; les
ensembles
$\;S=\{\,L(\omega)\,,\; \omega \in {\cal L}({\cal H})^*\,\}\;$ et
$\;\widehat S=\{\,\rho(\omega)\,,\; \omega \in {\cal L}({\cal
H})^*\,\}\;$ sont des sous-alg\`ebres involutives de $\;{\cal L}({\cal
H})\,$. Les applications $\;\delta:x\mt V(x\otimes 1)V^*\;$ de
$\;S\;$ dans
$\;S\otimes S\;$ et $\;\widehat\delta:x\mt V^*(1\otimes x)V\;$ de
$\;\widehat S\;$ dans $\;\widehat S\otimes \widehat S\;$
munissent
$\;S\;$ et $\;\widehat S\;$ de structures de $\,C^*$-alg\`ebres de
Hopf.

\subsec Remarquons que $\;L(\omega)\;$ (resp. $\;\rho(\omega)\,)$
ne d\'epend que de la restriction de $\;\omega \;$ \`a $\;\widehat
S\;$ (resp. $\;S\,$). Les alg\`ebres $\;S\;$ et $\widehat S\;$ sont en
dualit\'e {\it via}\/ la forme bilin\'eaire $\;\beta :(L(\omega
),\rho(\omega '))\mt (\omega\otimes \omega ') (V)=\omega (\rho
(\omega '))=\omega'(L(\omega ))\,$. Le produit (resp. coproduit) de
$\;S\;$ est dual pour $\;\beta\;$ du coproduit (resp. produit) de
$\;\widehat S\,$: pour $\;x,x'\in S\,,\;y,y'\in \widehat S\,$, on a
$$\beta (xx',y) =(\beta\otimes\beta)(x\otimes x',\widehat \delta
(y))\,,\qquad \beta (x,yy') =(\beta\otimes\beta)(\delta (x),y\otimes
y')\,.$$

Notons $\;\tau\;$ la trace normalis\'ee de $\;{\cal L}({\cal H})\,$.
Rappelons (\cf \eg [\ref\BSb], \S4) que les restrictions de $\;\tau
\;$ \`a
$\;S\;$ et $\;\widehat S\;$ sont des mesures de Haar: pour tout
$\;x\in S\;$ (resp $\;y\in \widehat S\,$) on a $\;(\id\otimes
\tau)\delta (x)=(
\tau\otimes\id)\delta (x)=\tau (x)1\;$ (resp. $\;(\id\otimes
\tau)\widehat \delta (y)=( \tau\otimes\id)\widehat \delta (y)=\tau
(y)1\,$).

Rappelons que $\;\xi\in{\cal H}\;$ est appel\'e {\it fixe}\/ (resp {\it
cofixe}) si, pour tout $\;\eta\in {\cal H}\;$ on a $\;V(\xi\otimes
\eta)=\xi\otimes \eta \;$ (resp. $\;V(\eta\otimes \xi)=\eta\otimes
\xi\,$)). Par [\ref\BSb] \S4, $\;L(\tau)\;$ (resp $\;\rho (\tau)\,$) est le
projecteur orthogonal sur l'espace des vecteurs cofixes (resp. fixes)
de $\;V\,$. Comme $\;(\tau\otimes \tau)(V)=\tau (L(\tau ))=\tau (\rho
(\tau ))\,$, les projecteurs $\;L(\tau )\;$ et  $\;\rho (\tau )\;$ ont
donc m\^eme dimension. Autrement dit, les dimensions de l'espace
des vecteurs fixes et des vecteurs cofixes co\"incident.

On appelle {\it multiplicit\'e}\/ de $\;V\;$ la dimension de l'espace
des vecteurs fixes. Rappelons (\cf [\ref\BSb]) qu'un unitaire
multiplicatif de multiplicit\'e $\;m\;$ est \'equivalent au produit
tensoriel d'un unitaire multiplicatif de multiplicit\'e $\;1\;$ par
$\;1\in \L(\C^m\otimes \C^m)\,$.

\subsec On suppose d\'esormais que la multiplicit\'e de $\;V\;$ est
$\;1\,$. Soient $\;e\;$ un vecteur fixe et $\;\hat e\in {\cal H}\;$ un
vecteur cofixe de norme 1. Soient $\;\varphi \;$ l'\'etat vectoriel
$\;\omega _{e,e}\;$ et $\;\widehat\varphi \;$ l'\'etat vectoriel
$\;\omega _{\hat e,\hat e}\,$. Rappelons que $\;\varphi \;$ et
$\;\tau\;$ co\"incident sur $\;S\;$ et que $\;\widehat\varphi\;$ et
$\;\tau\;$ co\"incident sur $\;\widehat S\,$. Comme
$\;L(\tau)=\theta_{\hat e,\hat e}\,$, on a $\;1/n=\tau (\theta_{\hat
e,\hat e})=\varphi (\theta_{\hat e,\hat e })=|\langle e,\hat e\rangle
|^2\,$. Quitte \`a multiplier $\;\hat e\;$ par un scalaire de module $\;1\;$
convenable, on supposera d\'esormais que $\; \langle e,\hat
e\rangle=n^{-1/2}\,$.

 Les applications $\;\kappa:L(\omega)\mt
L(\omega^*)^*=(\omega\otimes \id )(V^*) \;$ et $\;\widehat
\kappa:\rho(\omega)\mt \rho(\omega^*)^* = (\id\otimes\omega )(V^*)\;$
sont des antiautomorphismes involutifs de
$\;S\;$ et $\;\widehat S\,$ respectivement.

Comme $\;\kappa\;$ pr\'eserve l'involution de $\;S\,$, l'application
$\;xe\mt \kappa(x)e\;$ est un unitaire $\;U\in {\cal L}({\cal H})\,$.
On a $\;U^2=1\,$. De plus, $\;({\cal H},V,U)\;$ est un syst\`eme de Kac, au
sens de [\ref\BSb], \S 6. En particulier, en d\'esignant par
$\;\Sigma\in {\cal L}({\cal H}\otimes \H)\;$ la volte (\ie
l'op\'erateur tel que $\;\Sigma (\xi\otimes \eta)=\eta \otimes
\xi\,$), les unitaires $\;\widetilde V=(U\otimes 1)\Sigma V\Sigma
(U\otimes 1)\;$ et $\;\widehat V=\Sigma (U\otimes 1)V(U\otimes
1)\Sigma \;$ sont multiplicatifs. De plus,
$\;\Sigma \widehat VV\widetilde V(1\otimes U)=1\;$ et, pour
$\;x\in S\,$, on a $\;\delta (x)=\widehat V^*(1\otimes x)\widehat
V\;$ et pour
$\;x\in \widehat S\,$, on a $\;\widehat \delta (x)=\widetilde
V(x\otimes 1)\widetilde V^*\,$.

Remarquons que $\;Ue=e\,$; de plus, si $\;\omega \in {\cal L}({\cal
H})^*\;$ satisfait $\;\omega=\omega^*\;$ et $\;L(\omega)=L(\omega
)^*\,$, alors $\;U L(\omega)e=L(\omega)e\,$. Appliquant ceci \`a
$\;n^{1/2}\widehat\varphi\,$, on trouve $\;U\hat e=\hat e\,$.

Il s'ensuit que $\;\hat e\;$ est un vecteur fixe et $\;e\;$ un vecteur
cofixe pour
$\;\widetilde V\;$ et $\;\widehat V\,$.

\subsec On a $\;(\id\otimes \omega_{\hat e,e})(\Sigma
V)=(\id\otimes
\omega _{\hat e,e})(\Sigma \widehat VV\widetilde V(1\otimes
U))=n^{-1/2}\,$; de m\^eme, $\;(\omega_{e,\hat e}\otimes
\id)(\Sigma V)=(\omega_{e,\hat e}\otimes \id)(\Sigma \widehat
VV\widetilde V)=n^{-1/2}U\,$. Pour $\;\xi\in H\,$, on a donc
$\;L(\omega _{\hat e,\xi})e=(\id \otimes \omega_{\hat e,e})(\Sigma
V)\xi=n^{-1/2}\xi\;$ et
$\;\rho(\omega _{\xi,e})^*\hat e=(\id \otimes \omega_{\hat
e,e})(\Sigma V)^*\xi=n^{-1/2}\xi\,$. Comme $\;e\;$ et $\;\hat e\;$
sont s\'eparateurs pour $\;S\;$ et $\;\widehat S\;$ respectivement,
on en d\'eduit que, pour
$\;x\in S\,,\;y\in \widehat S\;$ on a $\;x=n^{1/2}L(\omega _{\hat
e,xe})\;$ et $\;y=n^{1/2}\rho (\omega _{y\hat
e,e})^*=n^{1/2}\widehat\kappa(\rho (\omega _{e,y\hat e}))\,$.

En particulier, $\;\beta (x,y)=\beta (n^{1/2}L(\omega _{\hat
e,xe}),y)=n^{1/2}\omega _{\hat e,xe}(y)=n^{1/2} \langle \hat e,
yxe\rangle\,$.

Remarquons de plus que, pour $\;y\in \widehat S\,,$ on a $\;\widehat
\kappa(y)=n^{1/2}\rho (\omega _{e,y\hat e})\,$, donc $\;\widehat
\kappa(y)\hat e=n^{1/2}\rho (\omega _{e,y\hat e})\hat
e=n^{1/2}(\omega_{e,\hat e}\otimes \id)(\Sigma V)y\hat e=Uy\hat e
\,$.

\ssec {La transformation de Fourier.} Pour $\;x\in S\,$, on pose
$\;\F(x)=n^{1/2}\rho (\omega _{e,xe })\,$; pour $\;y\in \widehat
S\,$, on pose $\;\widehat {\F}(y)=n^{1/2}L(\omega _{\hat e,y\hat
e})\,$. Par ce qui pr\'ec\`ede, pour $\;x\in S\;$ et $\;y\in \widehat
S\;$ on a $\;\F(x)\hat e=Uxe\;$ et $\;\widehat {\F}(y)e=y\hat e\,$. En
particulier, $\;\widehat {\F}\circ \F=\kappa\;$ et $\;\F\circ
\widehat {\F}=\widehat \kappa\,$. On en d\'eduit aussi que $\;\F
\circ \kappa = \F \circ \widehat {\F}\circ
\F=\widehat \kappa\circ \F\,$.

Nous caract\'erisons ci-dessous les \'el\'ements dont la
transform\'ee de Fourier est positive, puis ceux dont la
transform\'ee de Fourier est centrale.

\prop{a) Pour $\;x \in S\,$, les conditions suivantes sont
\'equivalentes:
\hfill\break
 (i) $\;\rho(\omega _{e,xe })\ge 0\,$. \quad
 (ii) La forme $\;\omega_{\hat e, xe}\;$ est positive sur $\;\widehat
S\,$. \quad (iii) Il existe une forme positive $\;\psi \in\L(\H)^*\;$
telle que $\;x=L(\psi)\,$.
\hfill\break
\indent L'ensemble des $\;x\in S\;$ tels que $\;\rho(\omega _{e,xe})
\ge 0\;$ est un c\^one convexe ferm\'e stable par le produit de
$\;S\;$ et invariant par $\;\kappa\,$.\hfill\break
 b) Pour $\;y \in \widehat S\,$, les conditions suivantes sont
\'equivalentes:
\hfill\break
 (i) $\;L(\omega _{\hat e,y\hat e })\ge 0\,$. \quad
 (ii) La forme $\;\omega_{e, y\hat e}\;$ est positive sur $\;S\,$.
\quad (iii) Il existe une forme positive $\;\psi \in\L(\H)^*\;$ telle
que
$\;y=\rho(\psi)\,$.
\hfill\break
\indent L'ensemble des $\;y\in \widehat S\;$ tels que $\;L(\omega
_{\hat e,y\hat e })\ge 0\;$ est un c\^one convexe ferm\'e stable par
le produit de $\;\widehat S\;$ et invariant par
$\;\widehat\kappa\,$.}

\pf a) Posons $\;y=\widehat \kappa(\rho(\omega _{e,xe }))\,$. On a
$\;\rho(\omega _{e,xe})\ge 0 \iff y\ge 0\iff$ la forme
$\,y\widehat\varphi\;$ est positive sur $\;\widehat S\,$. Or, $\;x
e=n^{1/2}y\hat e \,,$ d'o\`u l'\'equivalence (i)$\iff$(ii).

On a $\;x=nL(\omega_{\hat e, y\hat e})\;$ donc, pour une forme
$\;\psi\,$, on a $\;L(\psi)=x\;$ si et seulement si $\;\psi\;$
co\"incide avec $\;n^{1/2}\omega_{\hat e, xe}\;$ sur $\;\widehat
S\,$, d'o\`u (iii)$\Ra$(ii). La r\'eciproque r\'esulte de ce que toute
forme positive sur $\;\widehat S\;$ peut se prolonger en une forme
positive sur $\;\L(\H)\,$.

 Si $\;\omega \;$ et $\;\omega '\;$ sont des formes positives,
$\;\psi:x\mt (\omega \otimes \omega')(V^*(1\otimes x)V)\;$ est
positive et $\;L(\omega )L(\omega ')=L(\psi)\,$. En outre, pour
$\;x\in S\,$, si
$\;\F(x)\ge 0\;$ alors $\;\F(\kappa (x))= \widehat \kappa (\F(x))
\ge0\,$. L'ensemble des $\;x\in S\;$ tels que $\;\rho(\omega _{e,xe})
\ge 0\;$ est donc stable par le produit de $\;S\;$ et invariant par
$\;\kappa\,$. Les autres assertions de a) sont claires.

\medskip\noindent b) Se d\'emontre de mani\`ere analogue (ou en
rempla\c cant $\;V\;$ par $\;\widetilde V\,$).

\cqfd

Notons par la m\^eme lettre $\;\sigma \;$ la volte  de $\;S\otimes
S\;$ et celle de $\;\widehat S\otimes\widehat S\;$ (\ie
l'automorphisme donn\'e par $\;\sigma( x\otimes y)=y\otimes x\,$).

\prop {Soit $\;x\in S\;$ (resp. $\;\widehat S\,$). Alors
$\;\rho(\omega _{e,xe })\;$ (resp. $\;L(\omega _{\hat e,x\hat e })\,$)
est dans le centre de $\;\widehat S\;$ (resp $\;S\,$) si et seulement
si $\;\delta (x)=\sigma(\delta (x))\;$ (resp. $\;\widehat
\delta (x)=\sigma(\widehat \delta (x))\,$. L'ensemble des $\;x\in
S\;$ (resp. $\;\widehat S\,$) tels que $\;\rho(\omega _{e,xe})\;$
(resp. $\;L(\omega _{\hat e,x\hat e })\,$) est dans le centre de
$\;\widehat S\;$ (resp $\;S\,$) est une sous-alg\`ebre involutive de
$\;S\;$ (resp. $\;\widehat S\,$).}

\pf $\;\rho(\omega _{e,xe })\;$ (resp. $\;L(\omega _{\hat e,x\hat e
})\,$) est dans le centre de $\;\widehat S\;$ (resp $\;S\,$) si et
seulement si la forme $\;y\mt  \widehat \varphi (y\rho(\omega
_{e,xe }))\;$ (resp. $\;y\mt  \varphi (yL(\omega _{\hat e,x\hat e
}))\,$) est traciale sur $\;\widehat S\;$ (resp. $\;S\,$). Or,
$\;\widehat \varphi (y\rho(\omega _{e,xe }))=\langle\hat
e,y\rho(\omega _{e,xe })\hat e\rangle =n^{-1/2}\langle\hat
e,y\kappa (x)e\rangle=n^{-1}\beta (\kappa (x),y)\;$ (resp. $\;\varphi
(yL(\omega _{\hat e,x\hat e }))=\langle e,yL(\omega _{\hat e,x\hat e
})e\rangle =n^{-1/2}\langle e,yx\hat e\rangle =n^{-1}\overline{\beta
(y^*,x^*)}\,$). Cela signifie donc que pour $\;y,z\in \widehat S\;$
(resp $\;S\,$), on a $\;\beta (\kappa (x),yz)=\beta (\kappa (x),zy)\;$
(resp. $\;\beta (yz,x^*)=\beta (zy,x^*)\,$). Or $\;\delta
(\kappa(x))=(\kappa\otimes \kappa)\sigma \delta (x)\,$. La
proposition r\'esulte alors des formules de \sn 1.

\cqfd

\subsec On note $\;J:{\cal H}\rightarrow {\cal H}\;$ l'isom\'etrie
antilin\'eaire $\;xe\mt x^*e\;$ et $\;\widehat J:{\cal H}\rightarrow
{\cal H}\;$ l'isom\'etrie antilin\'eaire $\;L(\omega)e\mt
L(\omega^*)e\,$. Remarquons que $\;U=J\widehat J\,$. Les
applications $\;J\,,\;\widehat J\;$ et $\;U\;$ \'etant involutives, on
a aussi $\;U=\widehat JJ\,$. Pour $\;y\in
\widehat S\,$, on a $\;y\hat e=n^{1/2}L(\omega _{\hat e,y\hat
e})e=n^{1/2}L(y\widehat \varphi)e\;$ donc $\;\widehat Jy\hat
e=n^{1/2}L(\widehat\varphi y^*)e=n^{1/2}L(y^*\widehat \varphi
)e\,$, vu que
$\;\widehat \varphi\;$ est une trace sur $\;\widehat S\,$. Donc
$\;\widehat Jy\hat e=y^*\hat e\,$. Enfin, $\;J\rho (\omega )\hat
e=\widehat JU\rho (\omega )\hat e=\widehat J\rho (\omega
^*)^*\hat e=\rho (\omega ^*)\hat e\,$.

Remarquons que, pour $\;x\in S\,$, $\;JxJ\in S'\,$, $\;UxU\in S'\;$ et
$\;\widehat Jx\widehat J=\kappa(x^*)\in S\,$. De m\^eme, pour
$\;y\in
\widehat S\,$, $\;\widehat Jy\widehat J\in \widehat S'\,$,
$\;UyU\in
\widehat S'\;$ et $\;JyJ=\widehat \kappa(y^*)\in \widehat S\,$.
Donc pour
$\;\omega\in \L(\H)^*\,$, la restriction \`a $\;\widehat S\;$ (resp.
$\;S
\,$) de la forme $\;x\mt \overline{\omega (JxJ)}\;$ (resp. $\;x\mt
\overline{\omega (\widehat Jx\widehat J)}\,$) ne d\'epend que de la
restriction de $\;\omega \;$ \`a $\;\widehat S\;$ (resp. $\;S\,$). Pour
$\;\xi\in \H\,$, on a $\;L(\omega_{\hat e,\xi})^*e=JL(\omega_{\hat
e,\xi})e=n^{-1/2}J\xi=L(\omega_{\hat e,J\xi})e\,$, donc $\;L(\omega
_{\hat e,\xi})^*=L(\omega_{\hat e,J\xi})\,$. De m\^eme, $\;\rho
(\omega _{e,\xi})^*=\rho (\omega _{e,\widehat J\xi})\,$. Comme la
restriction de toute forme \`a $\;\widehat S\;$ (resp. $\;S\,$)
co\"incide avec un
$\;\omega _{\hat e,\xi}\;$ (resp. un $\;\omega _{e,\xi}\,$), on en
d\'eduit que, pour tout $\;\xi,\eta\in\H\,$, on a:
$\;L(\omega_{\xi,\eta })^*=L(\omega _{J\xi,J\eta })\;$ et
$\;\rho(\omega_{\xi,\eta })^*=\rho (\omega _{\widehat
J\xi,\widehat J\eta })\,$.

\prop{a) Pour tout $\;x\in S\,,\;y\in \widehat S\;$ on a: $\;\tau
(xy)=\tau (x)\tau (y)= \varphi (x)\widehat\varphi(y)\,$. \hfill\break
 b) On a $\;(\id\otimes \tau \otimes \id)(V_{12}V_{23})=(\id\otimes
\tau \otimes \id)(V_{23}V_{12})=\theta _{e,e}\otimes \theta _{\hat
e,
\hat e}\,$. \hfill\break
 c) Pour tout $\;x\in S\,,\;y\in \widehat S\;$ on a: $\;L( x\tau)=\tau
(x)
\theta _{\hat e,\hat e}\;$ et $\;\rho (y\tau )=\tau (y)\theta
_{e,e}\,$.}

\pf a) On a $\;\widetilde V(xy\otimes 1)\widetilde V^*=(x\otimes
1)\widehat \delta (y)\,$. Donc  $$\eqalign{\tau (xy)&=(\tau \otimes
\tau)(xy\otimes 1)\cr &=(\tau \otimes \tau)(\widetilde V(xy\otimes
1)\widetilde V^*) \qquad\hbox{vu que $\;\tau \otimes \tau \;$ est
une trace}\cr &=(\tau \otimes \tau)((x\otimes 1)\widehat \delta
(y))\cr &=\tau(x(\id\otimes \tau)\widehat \delta (y))\cr
&=\tau(x)\tau (y)\,.}$$

\noindent b) Pour $\;\omega,\omega'\in {\cal L}({\cal H})^*\,$, on a
$\;(\omega
\otimes \omega ')(\id\otimes \tau \otimes
\id)(V_{12}V_{23})=\tau(\omega
\otimes \id\otimes \omega ')(V_{12}V_{23})=\tau (L(\omega )\rho
(\omega '))=\tau(L(\omega ))\tau (\rho (\omega '))\;$ par a). Or
$\;\tau(L(\omega ))=\omega (\rho(\tau))\;$ et $\;\tau(\rho(\omega
'))=\omega' (L(\tau))\,$. Donc $\;(\omega \otimes \omega
')(\id\otimes \tau
\otimes \id)(V_{12}V_{23})=(\omega \otimes \omega ')(\theta
_{e,e}\otimes
\theta _{\hat e, \hat e})\,$.

\noindent c) Pour $\;\omega\in {\cal L}({\cal H})^*\,$, on a  $\;\omega
(L(\tau x))=\tau (x\rho (\omega ))=\tau(x)\omega(L(\tau ))\,$; de
m\^eme,
$\;\omega (\rho(\tau y))=\tau (yL (\omega
))=\tau(y)\omega(\rho(\tau ))
\,$.

\cqfd

\cor{Si $\;p\in S\;$ et $\;q\in \widehat S\;$ sont des projecteurs qui
commutent, on a $\;\dim(p\H)\dim(q\H)=n\dim(p\H\cap
q\H)=n^2\|pe\|^2\|q\hat e\|^2\,$.}

\pf Cela r\'esulte du calcul de $\;\tau(pq)\;$ \`a l'aide de la prop. \sn
7.a\ref).

\cqfd

Il est facile de donner une description de l'esp\'erance
conditionnelle de $\;\L(\H)\;$ sur $\;S\,$:

\prop{Pour tout $\;T\in \L(\H)\;$ on a $\;E_S(T)=(\id\otimes \omega
_{\hat e,\hat e})(\widehat V^*(1\otimes T)\widehat V)\,$.}

\pf Il suffit de v\'erifier cette formule pour $\;T=ab\;$ avec $\;a\in
S\;$ et $\;b\in \widehat S\,$. On a $\;(\id\otimes \omega _{\hat
e,\hat e})(\widehat V^*(1\otimes T)\widehat V)= (\id\otimes
\omega _{\hat e,b\hat e})(\widehat V^*(1\otimes a)\widehat V)=
(\id\otimes \omega _{\hat e,b\hat e})(V(a\otimes 1)V^*)=\tau
(b)a\,$.

\cqfd

On peut donner d'autres \'ecritures pour $\;E_S\;$ par exemple
$\;E_S(T)=(\omega _{\hat e,\hat e}\otimes \id )(V(T\otimes
1)V^*)=(\id\otimes \tau)(\widetilde V(T\otimes 1)\widetilde
V^*)\,$, ainsi que des formules analogues pour $\;E_{\widehat S}\,$.

On peut remarquer que l'unitaire multiplicatif $\;V\;$ (de
multiplicit\'e $\;1\,$) est enti\`erement d\'etermin\'e par les alg\`ebres
$\;S\;$ et $\;\widehat S\;$ associ\'ees et les espaces de ses vecteurs
fixes et cofixes: par \sn3, pour $\;\omega \in \L(\H)^*\,$, l'op\'erateur
$\;\rho (\omega )\;$ est l'unique \'el\'ement $\;y\;$ de $\;\widehat
S\;$ tel que, pour tout $\;x\in S\;$ on ait $\;\omega
(x)=n^{1/2}\langle \hat e,yxe\rangle \,$. On peut donner une
description plus explicite de $\;V\;$ en termes de $\;S\,$,
$\;\widehat S\,$, $\;e\;$ et $\;\hat e\,$:

\prop{Pour tout $\;a,x\in S\;$ et tout $\;b,y\in \widehat S\,$, on a
$\;\langle xe\otimes y\hat e,V(ae\otimes b\hat e)\rangle=\tau
(x^*y^*ab)\,$.}

\pf On a $\;\langle xe\otimes y\hat e,V(ae\otimes b\hat e)\rangle=
\langle xe\otimes y\hat e,\delta (a)(e\otimes b\hat e)\rangle=
\langle xe,(\id \otimes \omega _{y\hat e,b\hat e})(\delta
(a))e\rangle\,$. Or $\;(\id \otimes \omega _{y\hat e,b\hat e})(\delta
(a))= (\id \otimes \omega _{\hat e,\hat e})(\widehat V^*(1\otimes
y^*ab)\widehat V)=E_S(y^*ab)\,$; donc $\;\langle xe,(\id \otimes
\omega _{y\hat e,b\hat e})(\delta (a))e\rangle=\varphi
(x^*E_S(y^*ab))=\tau (x^*y^*ab)\,$.

\cqfd

\cor{On a $\;(J\otimes\widehat J)V(J\otimes\widehat J)=V^*\,$.}

\pf Si $\;a,x\in S\;$ et $\;b,y\in \widehat S\,$, on a
$\; \langle (J\otimes\widehat J)V(J\otimes\widehat J) (xe\otimes
y\hat e),  ae\otimes b\hat e\rangle=\langle (J\otimes\widehat
J)(ae\otimes b\hat e),V(J\otimes\widehat J) (xe\otimes y\hat
e)\rangle =\langle a^*e\otimes b^*\hat e,  V(x^*e\otimes y^*\hat
e)\rangle =\tau (abx^*y^*)=
\langle xe\otimes y\hat e,V(ae\otimes b\hat e)\rangle\,$.

\cqfd

\sec {Une caract\'erisation des alg\`ebres de Hopf en dualit\'e.}

Soit $\;\H\;$ un espace hilbertien de dimension finie (non nulle)
$\;n\,$. On note $\;\tau \;$ la trace normalis\'ee de $\;\L(\H)\,$. Soient
$\;A,B\subset \L(\H)\;$ deux sous-alg\`ebres involutives unitales. Nous
donnons ici des conditions pour que $\;A\;$ et $\;B\;$ soient des alg\`ebres
de Hopf en dualit\'e, \ie les alg\`ebres $\;S\;$ et $\;\widehat S\;$
associ\'ees \`a un unitaire multiplicatif. La premi\`ere condition que nous
imposons est
une condition de \og carr\'e commutatif\fg (\cf [\Po]).

\lem{Soient $\;A,B\subset \L(\H)\;$ deux sous-alg\`ebres involutives
unitales telles que, pour tout  $\;(a,b)\in A\times B\,$, on ait
$\;\tau (ab)=\tau (a)\tau (b)\,$.\hfill\break
 a) Soit $\;\hat f\in \H\;$ de norme $\;1\,$. Si $\;\theta _{\hat f,\hat f}\in
A\;$ alors $\;\hat f\;$ est un vecteur trace pour $\;B\,$. \hfill\break
 b) Si de plus $\;\dim B=n\;$ alors $\;\theta _{\hat f,\hat f}\;$ est central
dans $\;A\,$.}

\pf a)) Pour $\;b\in B\,$, on a $\;\langle \hat f,b\hat f\rangle=n\tau
(\theta_{\hat f,\hat f}b)=n\tau (\theta_{\hat f,\hat f})\tau (b)=\tau (b) \,$,
d'o\`u a).

\noindent b) Par a), $\;\hat f\;$ est un vecteur s\'eparateur pour
$\;B\,$; comme $\;\dim B\ge \dim\,\H\,$, on en d\'eduit que c'est un
vecteur totalisateur pour $\;B\,$.

Soit $\;a\in A\,$. Pour tout $\;b\in
B\,$, on a $\;\langle b^*\hat f,a\hat f\rangle =n\tau
(a\theta _{\hat f,\hat f}b)=n\tau (a\theta _{\hat f,\hat f})\tau (b)=\langle
b^*\hat f,\hat f\rangle\langle \hat f,a\hat f\rangle\;$ donc $\;a\hat
f=\langle \hat f,a\hat f\rangle \hat f\,$ (puisque $\;\hat f\;$ est un
totalisateur pour $\;B\,$), donc $\;a\theta_{\hat f,\hat f}=\langle \hat
f,a\hat f\rangle \theta_{\hat f,\hat f}\,$. Appliquant cela \`a $\;a^*\,$, on
trouve $\;a\theta_{\hat f,\hat f}=\theta _{\hat f,\hat f}a\,$.

\cqfd

Soient $\;A,B\subset \L(\H)\;$ deux sous-alg\`ebres involutives
unitales de dimension $\;n\;$ telles que, pour tout  $\;(a,b)\in
A\times B\,$, on ait $\;\tau (ab)=\tau (a)\tau (b)\;$ et $\;\hat f,f\in
\H\;$ de norme $\;1\;$ tels que $\;\theta_{\hat f,\hat f}\in A\;$ et
$\;\theta_{f,f}\in B\,$. Remarquons que, comme $\;\hat f\;$ est un
vecteur trace pour $\;B\,$, on a $\;|\langle \hat f,f\rangle
|=n^{-1/2}\,$. On supposera d\'esormais que $\;\langle \hat f,f\rangle
=n^{-1/2}\;$ (quitte \`a multiplier $\;\hat f\;$ par un scalaire de
module $\;1\,$).

On note $\;\Lambda\;$ l'injection canonique de $\;\L(\H)\;$ dans
$\;L^2(\L(\H),\tau )\,$. D\'efinissons $\;Z,Z':\H\otimes \H\ra
L^2(\L(\H),\tau )\;$ par les formules $\;Z(af\otimes b\hat f)=\Lambda
(ab)\;$ et $\;Z'(af\otimes b\hat f)=\Lambda (ba)\,$. Il est clair que ce
sont des isom\'etries, donc des unitaires (par \'egalit\'e des
dimensions). On pose enfin $\;W=(Z')^*Z\,$.

\lem{a) Pour tout $\;a,x\in A\;$ et tout $\;b,y\in B\,$, on a
$\;\langle xf\otimes y\hat f,W(af \otimes b\hat f)\rangle =\tau
(x^*y^*ab)\,$.\hfill\break
 b) Pour tout $\;b\in B\;$ et tout $\;a\in A\,$, on a $\;n^{1/2}(\omega
_{\hat f,af} \otimes \id)(W)=a\;$ et $\;n^{1/2}(\id \otimes \omega
_{b^*\hat f,f})(W)=b\,$.\hfill\break
 c) Les propri\'et\'es suivantes sont \'equivalentes: (i) $\;W\in
B\otimes \L(\H)\,$; (ii) $\;W\in \L(\H)\otimes A\,$; (iii) $\;W\in
B\otimes A\,$.}

\pf a) En effet, $\;\langle xf\otimes y\hat f,W(af \otimes b\hat f)\rangle=
\langle Z'(xf\otimes y\hat f),Z(af \otimes b\hat f)\rangle=\langle
\Lambda(yx),\Lambda (ab)\rangle \,$.

\noindent b) Pour $\;a\in A\;$ et $\;b,y\in B\,$, on a $\;\langle
y\hat f,n^{1/2}(\omega _{\hat f,af} \otimes \id)(W)b\hat
f\rangle =n^{1/2}\langle
\hat f\otimes  y\hat f,W(af\otimes b\hat f)\rangle =n\langle \theta_{\hat
f,\hat f}f\otimes  y\hat f,W(af\otimes b\hat f)\rangle=n\tau
(\theta _{\hat f,\hat f}y^*ab)=\langle y\hat f,ab\hat f\rangle \,$. Comme
$\;\hat f\;$ est totalisateur pour $\;B\,$, on en d\'eduit la premi\`ere
assertion.

En \'echangeant les r\^oles de $\;A\;$ et $\;B\,$, on remplace
$\;W\;$ par $\;\Sigma W^*\Sigma\,$, d'o\`u la deuxi\`eme assertion.

\noindent c) Comme $\;f\;$ est totalisateur pour $\;A\;$ et $\;\hat f\;$
s\'eparateur pour $\;B\,$, tout \'el\'ement de $\;B^*\;$ est la restriction
d'une forme $\;\omega _{\hat f,af}\,$. Si $\;W\in B\otimes \L(\H)\,$, pour
tout  $\;\omega \in B^*\,$, on a $\;(\omega \otimes \id)(W)\in A\;$ par b),
donc $\;W\in B\otimes A\,$. On a montr\'e (i)$\iff$(iii). Echangeant les
r\^oles de $\;A\;$ et $\;B\,$, on trouve que (ii)$\iff$(iii).

\cqfd

Remarquons que sous les hypoth\`eses du lemme \sn2.c), on a
$\;B=\{\,(\id\otimes \omega)(W)\,,\;\omega\in\L(\H)^*\,\}\;$ et $\;A=\{\,(
\omega \otimes \id)(W)\,,\;\omega\in\L(\H)^*\,\}\,$.

Notons $\;J_A\;$ et $\;J_B\;$ les involutions de $\;\H\;$d\'efinies
par $\;J_A(af)=a^*f\;$ et $\;J_B(b\hat f)=b^*\hat f\;$ (pour $\;a\in A\;$ et
$\;b\in B\,$).

\lem{a) On a $\;(J_A\otimes J_B)W(J_A\otimes J_B)=W^*\,$.
\hfill\break
 b) Si $\;W\in \L(\H)\otimes A\,$, alors $\;J_ABJ_A=B\;$ et
$\;J_BAJ_B=A\,$.\hfill\break
 c) Dans ce cas, $\;J_AJ_B=J_BJ_A\,$.}

\pf a) Pour $\;a,x\in A\;$ et $\;b,y\in B\,$, on a $\;\langle
(J_A\otimes J_B)W(J_A\otimes J_B)(xf\otimes y\hat f),af\otimes b\hat f
\rangle =\langle (J_A\otimes J_B)(af\otimes b\hat f),W(J_A\otimes
J_B)(xf\otimes y\hat f)\rangle =\langle a^*f\otimes b^*\hat f,W(x^*f
\otimes y^*\hat f)\rangle =\tau (x^*y^*ab)=\langle xf\otimes y\hat
f,W(af\otimes b\hat f)\rangle\,$, d'o\`u a).

\noindent b) Pour $\;b\in B\,$, on a $\;J_AbJ_A=n^{1/2}J_A((\id\otimes
\omega _{b^*\hat f,f})(W))J_A=n^{1/2}(\id\otimes (J_B\omega
_{b^*\hat f,f}J_B))(W^*)\in B\,$, d'o\`u la premi\`ere assertion. La
deuxi\`eme assertion se d\'emontre de fa\c con analogue.

\noindent c) On a $\;J_A(\hat f)=n^{1/2}J_A\theta _{\hat f,\hat f}f=\hat
f\,$, donc $\;J_BJ_Ab\hat f=J_B(J_AbJ_A)\hat f=(J_AbJ_A)^*\hat
f=J_AJ_Bb\hat f\,$.

\cqfd

On suppose d\'esormais que $\;W\in \L(\H)\otimes A\,$. On pose
$\;U=J_AJ_B\;$ et $\;\widehat W=\Sigma (U\otimes 1)W(U\otimes
1)\Sigma \,$.

\lem{Pour $\;x,a\in A\;$ et $\;b',y'\in B'\,$, on a $\;\langle
y'\hat f\otimes xf,\widehat W(b'\hat f \otimes af)\rangle =\tau
(y'{}^*x^*b'a)\,$.}

\pf On a $\;\langle y'\hat f\otimes xf,\widehat W(b'\hat f \otimes
af)\rangle =\langle Uxf\otimes y'\hat f,W(Uaf \otimes b'\hat
f)\rangle \,$. Posons $\;a_1=J_Ba^*J_B\,,
\;x_1=J_Bx^*J_B\,,\;b=J_B(b')^*J_B\;$ et $\;y=J_B(y')^*J_B\,$. On a
$\;a_1f=Uaf\,,\;x_1f=Uxf\,,\; b\hat f=b'\hat f\;$ et $\;y\hat f=y'\hat f\,$.
Donc $\;\langle y'\hat f\otimes xf,\widehat W(b'\hat f \otimes af)\rangle
=\langle x_1f\otimes y\hat f,W(a_1f \otimes b\hat f)\rangle=\tau
(x_1^*y^*a_1b)=\tau (J_B(b'ay'{}^*x^*)^*J_B)=\tau (b'a y'{}^*x^*)\,$.

\cqfd

\lem{a) Pour tout $\;a\in A\;$ on a $\;Z(a\otimes 1)Z^*=\pi _\tau
(a)\,$.\hfill\break
 b) Pour tout $\;a_1,a_2\in A\;$ et $\;b_1,b_2\in B\,$, on a $\;\langle
a_1f\otimes b_1\hat f,W(\theta _{\hat f,\hat f}\otimes 1)W^*(a_2f\otimes
b_2\hat f)\rangle= n^{-1}\langle \hat f,b_2a_2a_1^*b_1^*\hat f\rangle\,$.
\hfill\break
 c) Pour tout $\;c_1,c_2\in A\;$ et $\;d_1,d_2\in B'\,$, on a
$\;\langle d_1\hat f\otimes c_1f,\widehat W^*(1\otimes \theta
_{\hat f,\hat f})\widehat W(d_2\hat f\otimes c_2f)\rangle =n^{-1}\langle
\hat f,d_2c_2c_1^*d_1^*\hat f\rangle\,$.
\hfill\break
 d) On a $\;W(\theta_{\hat f,\hat f}\otimes 1)W^*=\widehat W^*(1\otimes
\theta _{\hat f,\hat f})\widehat W\,$.}

\pf a) est clair.

\noindent b) Comme $\;Z'(a_if\otimes b_i\hat f)=\Lambda_\tau
(b_ia_i)\,$, on a
$$\eqalign{\langle a_1f\otimes
b_1\hat f,W(\theta _{\hat f,\hat f}\otimes 1)W^*(a_2f\otimes b_2\hat
f)\rangle &=\langle(\Lambda_\tau (b_1a_1),Z(\theta _{\hat f,\hat f}\otimes 1)Z^*
\Lambda_\tau (b_2a_2)\rangle \cr & =\tau (a_1^*b_1^*\theta _{\hat f,\hat
f}b_2a_2)\qquad \hbox{ par a).}\cr& =n^{-1}\langle \hat
f,b_2a_2a_1^*b_1^*\hat f\rangle }$$
c) Remarquons que, comme $\;A=J_BAJ_B\;$ et $\;B'=J_BBJ_B\,$, pour tout
$\;a\in A\;$ et tout $\;b'\in B'\;$ on a $\;\tau (ab')=\tau (a)\tau (b')\,$. De
plus, $\;\theta _{f,f}\in B'\;$ et $\;\theta _{\hat f,\hat f}\in A\,$.
Si on remplace le couple $\;(A,B)\;$ par $\;(B',A)\;$ on remplace $\;W\;$
par $\;\widehat W\;$ (lemme \sn4\ref); donc si on remplace le couple
$\;(A,B)\;$ par $\;(A,B')\;$ on remplace $\;W\;$ par $\;\Sigma\widehat
W^*\Sigma \,$; donc c) d\'ecoule de b).

\noindent d) Soient $\;a_1,a_2\in A\;$ et
$\;b_1,b_2\in B\,$; notons $\;c_1,c_2\in A\;$ et $\;d_1,d_2\in B'\;$ les
\'el\'ements tels que $\;a_if=d_i\hat f\;$ et $\;b_i\hat f=c_if\;$
($\,i=1,2\,$). Remarquons que  $\;Uc_i^*f=J_Bc_if= b_i^*\hat f\;$ et
$\;Ua_i^*f=d_i^*\hat f\,$. On a: $$\eqalign{\langle a_1f\otimes b_1\hat
f,W(\theta _{\hat f,\hat f}\otimes 1)W^*(a_2f\otimes b_2\hat f)\rangle
&=n^{-1}\langle \hat f,b_2a_2a_1^*b_1^*\hat f\rangle \cr& =n^{-1}\langle
a_2^*Uc_2^*f,a_1^*Uc_1^*f\rangle \cr& =n^{-1} \langle
a_2^*f,U(c_2c_1^*)Ua_1^*f\rangle\quad \hbox{car $\;Uc_iU\in A'\,$.}\cr&
=n^{-1}\langle d_2^*\hat f,(c_2c_1^*)d_1^*\hat f\rangle\cr&= \langle
d_1\hat f\otimes c_1f,\widehat W^*(1\otimes \theta _{\hat f,\hat
f})\widehat W(d_2\hat f\otimes c_2f)\rangle}$$ d'o\`u le r\'esultat.

\cqfd

\lem{Notons $\;\pi:\L(\H)\ra \L(\H\otimes \H)\;$ la repr\'esentation
$\;x\mt (Z')^*\pi_\tau (x)Z'\,$. \hfill\break
 a) Pour $\;a\in A\;$ et $\;b\in B\,$,
on a $\;\pi (a)=W(a\otimes 1)W^*\;$ et $\;\pi (b)=(1\otimes
b)\,$.\hfill\break
 b) Pour tout $\;x\in \L(\H)\;$ on a $\;\pi(x)=\widehat W^*(1\otimes x)
\widehat W\,$.}

\pf a) La deuxi\`eme assertion est claire; la premi\`ere r\'esulte du
lemme \sn 5.a).

\noindent b) Comme $\;\widehat W\in A\otimes B'\,$, cette \'egalit\'e est
vraie pour $\;x\in B\;$ par a). Par le lemme \sn5.d), cela est aussi vrai
pour $\;x=\theta _{\hat f,\hat f}\,$. Or $\;\L(\H)\;$ est engendr\'e par
$\;B\;$ et $\;\theta _{\hat f,\hat f}\,$, d'o\`u le r\'esultat.

\cqfd

\prop{Pour tout $\;a\in A\;$ on a $\;W(a\otimes 1)W^*\in A\otimes A\,$;
pour tout $\;b\in B\;$ on a $\;W^*(1\otimes b)W\in B\otimes B\,$.}

\pf Par hypoth\`ese $\;W(a\otimes 1)W^*\in \L(\H)\otimes A\,$; par le
lemme \sn6, on a $\;W(a\otimes 1)W^*=\widehat W^*(1\otimes
a)\widehat W\in A\otimes \L(\H)\,$, d'o\`u la premi\`ere assertion. Quand
on \'echange les r\^oles de $\;A\;$ et $\;B\,$, $\;W\;$ est remplac\'e par
$\;\Sigma W^*\Sigma\,$; la deuxi\`eme assertion r\'esulte donc de la
premi\`ere.

\cqfd

\th{Soient $\;A,B\subset \L(\H)\;$ deux sous-alg\`ebres
involutives unitales de dimension $\;n\;$ et $\;\hat f,f\in \H\;$
de norme $\;1\,$. On suppose que:\hfill\break
 a) Pour tout  $\;(a,b)\in A\times B\,$, on a $\;\tau (ab)=\tau
(a)\tau (b)\,$.\hfill\break
 b) $\;\theta _{\hat f,\hat f}\in A\;$ et $\;\theta _{f,f}\in
B\,$.\hfill\break
 c) L'unitaire $\;W\in \L(\H\otimes \H)\;$ d\'efini par $\;\langle
xf\otimes y\hat f,W(af \otimes b\hat f)\rangle =\tau (x^*y^*ab)\;$ pour tout
$\;a,x\in A\;$ et tout $\;b,y\in B\,$, est contenu dans
$\;B\otimes \L(H)\,$.\hfill\break
Alors $\;W\;$ est multiplicatif; les alg\`ebres \og$\,S\,$\fg\ et
\og$\,\widehat S\,$\fg\ associ\'ees sont $\;A\;$ et $\;B\,$; le
vecteur $\;f\;$ est fixe et le vecteur $\;\hat f\;$ cofixe pour $\;W\,$.}

\pf On a d\'ej\`a montr\'e que $\;\{\,(\omega \otimes \id)(W)\,,\; \omega
\in \L(\H)^*\,\}=A\;$ et $\;\{\,(\id\otimes \omega)(W)\,,\;
\omega \in \L(\H)^*\,\}=B\,$. Il est clair que pour tout $\;\xi \in \H\,$, on a
$\;W(f\otimes \xi )=f\otimes \xi \;$ et $\;W(\xi \otimes \hat f)=\xi
\otimes \hat f\,$. La seule chose \`a montrer est que $\;W\;$ est
multiplicatif.

Par la prop. \sn7, on a $\;W_{23}W_{12}W_{23}^*\in B\otimes A\otimes
A\;$ et $\;W_{12}^*W_{23}W_{12}\in B\otimes B\otimes A\,$; donc
$\;W_{12}^*W_{23}W_{12}W_{23}^*\in B\otimes \C\otimes A\,$.

Pour $\;x\in A\;$ et $\;\eta \in \H\,$, comme $\;W(\hat f\otimes \hat
f)=\hat f\otimes \hat f\;$ et $\;W^*(f\otimes \eta)=f\otimes \eta\,$, on a
$$\eqalign{\langle
\hat f\otimes \hat f\otimes f, W_{12}^*W_{23}W_{12}W_{23}^*(xf\otimes
f\otimes \eta)\rangle &= \langle \hat f\otimes \hat f\otimes f,
W_{23}W_{12}(xf
\otimes f\otimes \eta)\rangle \cr&=\langle \hat f\otimes f, W((\omega
_{\hat f,xf}\otimes \id)(W)\otimes 1)(f\otimes \eta)\rangle\cr&=n^{-1/2}
\langle \hat f\otimes f, W(xf\otimes \eta)\rangle\quad
\hbox {lemme \sn2.b)}\cr&=\langle
\hat f\otimes \hat f\otimes f, W_{13}(xf\otimes f\otimes
\eta)\rangle}$$ d'o\`u l'on d\'eduit que  $\;
W_{12}^*W_{23}W_{12}W_{23}^*=W_{13} \,,$ vu que
$\;\hat f\otimes f\;$ est s\'eparateur pour $\;B\otimes A\,$.

\cqfd

Remarquons que, comme la dimension de $\;A\;$ est \'egale \`a celle de
$\;\H\,$, l'unitaire multiplicatif $\;W\;$ est irr\'eductible. En fait,
comme $\;U=J_AJ_B\,$, $\;(\H,W,U)\;$ est un syst\`eme de Kac au sens de
[\BSb].

\sec{Pr\'e-sous-groupes.}

Soient $\;{\cal H}\;$ un espace hilbertien de dimension finie $\;n\ne
0\;$ et $\;V\in {\cal L}({\cal H}\otimes {\cal H})\;$ un unitaire
multiplicatif de multiplicit\'e $\;1\;$ sur $\;{\cal H}\,$. On choisit un
vecteur fixe $\;e\;$ de module $\;1\;$ et on note $\;\hat e\;$ le vecteur
cofixe tel que $\;\langle e,\hat e\rangle=n^{-1/2}\,$.

Pour $\;\xi \in{\cal H}\,$, posons $\;H_\xi =\{\,\eta \in {\cal
H}\,,\;V(\eta \otimes \xi )=\eta \otimes \xi \,\}\;$ et $\;H^\xi =\{\,\eta
\in {\cal H}\,,\;V(\xi \otimes \eta )=\xi \otimes \eta \,\}\,$.

\lem{Soit $\;\xi \in{\cal H}\;$ tel que $\;\|\xi\|=1\,$. Notons $\;\psi \;$
l'\'etat vectoriel $\;\omega _{\xi ,\xi }\,$.\hfill\break
 a) On a $\;H^\xi =\{\,\eta \in {\cal H}\,,\;L (\psi )\eta =\eta\}\;$ et
$\;H_\xi =\{\,\eta \in {\cal H}\,,\;\rho (\psi )\eta =\eta
\}\,$.\hfill\break
  b) On a $\;V(H_\xi \otimes H_\xi )\subset {\cal H}\otimes H_\xi \;$ et
$\;V^*(H^\xi \otimes H^\xi )\subset H^\xi \otimes {\cal H}\,$.}

\pf a) Comme $\;\|L(\psi)\|\le 1\,$, on a: $\;L(\psi )\eta =\eta \iff
\langle
\eta ,L(\psi )\eta \rangle =\|\eta \|^2 \iff \langle \xi \otimes \eta ,V(\xi
\otimes \eta )\rangle =\|\eta \|^2 \iff V(\xi \otimes \eta )=\xi \otimes
\eta \,$. La m\^eme d\'emonstration s'applique pour $\;\rho (\psi )\,$.

\noindent b) Soient $\;\eta ,\zeta \in H_\xi \,$. Alors $\;\eta \otimes
\zeta \otimes \xi\;$ est invariant par $\;V_{13}\;$ et $\;V_{23}\;$ donc
par $\;V_{13}V_{23}=V_{12}^*V_{23}V_{12}\,$. Donc $\;V(\eta
\otimes \zeta )\otimes \xi\;$ est invariant par $\;V_{23}\;$ \ie
$\;V(\eta \otimes \zeta )\in {\cal H}\otimes H_\xi \,$.

La deuxi\`eme assertion se d\'eduit de la premi\`ere en rempla\c cant
$\;V\;$ par $\;\Sigma V^*\Sigma \,$.

\cqfd

\lem{a) Soient $\;\xi_1,\xi_2,\eta_1,\eta_2\in \H\;$ tels que
$\;V(\xi_1 \otimes \eta_1)=\xi_1\otimes \eta_1\,$ et
$\;V(\xi_2\otimes
\eta_2)=\xi_2\otimes \eta_2\,$. On a $\;L(\omega _{\xi_1, \xi _2})
L(\omega _{\eta_1, \eta _2}) =\langle \xi_1, \xi _2 \rangle L(\omega
_{\eta_1, \eta _2})\;$ et $\;\rho (\omega _{\xi_1, \xi _2}) \rho (\omega
_{\eta_1, \eta _2}) =\langle \eta _1, \eta _2 \rangle \rho (\omega
_{\xi_1,
\xi _2})\,$.\hfill\break b) Soient $\;f,g\in \H\;$ tels que $\;V(f\otimes
g)=f\otimes g\,$. Pour tout $\;\xi ,\eta \in{\cal H}\;$ les op\'erateurs
$\;L(\omega _{\xi ,f})\;$ et
$\;\rho (\omega _{g,\eta })^*\;$ commutent.}

\pf a) On a $\;L(\omega _{\xi_1, \xi _2}) L(\omega _{\eta_1, \eta _2})
=L(\omega)\,$, o\`u $\;\omega(x)=(\omega _{\xi_1, \xi _2}\otimes
\omega _{\eta_1, \eta _2})(V^*(1\otimes x)V)=\langle V(\xi_1 \otimes
\eta_1),(1\otimes x)V(\xi_2 \otimes \eta_2)\rangle = \langle \xi_1,
\xi _2 \rangle  \omega _{\eta_1, \eta _2}(x) \,$. De m\^eme, $\;\rho
(\omega _{\xi_1, \xi _2}) \rho (\omega _{\eta_1, \eta _2}) =\rho
(\omega)\,$, o\`u $\;\omega(x)=(\omega _{\xi_1, \xi _2}\otimes
\omega _{\eta_1, \eta _2})(V(x\otimes 1)V^*)=\langle V^*(\xi_1
\otimes \eta_1),(x\otimes 1)V^*(\xi_2 \otimes \eta_2)\rangle =
\langle \eta_1, \eta _2 \rangle  \omega _{\xi_1, \xi _2}(x) \,$.

\noindent b) On a $$\eqalign{L(\omega _{\xi ,f})\rho (\omega _{g,\eta
})^* &= (\omega _{\xi ,f}\otimes \id \otimes \omega _{\eta
,g})(V_{12}V_{23}^*) \cr &=(\omega _{\xi ,f}\otimes \id \otimes \omega
_{\eta ,g})(V_{23}^*V_{12}V_{13})\cr &=(\omega _{\xi ,f}\otimes \id
\otimes \omega _{\eta ,g})(V_{23}^*V_{12})\cr}$$ vu que  $\;V(f\otimes
g)=f\otimes g\,$. Or $\;(\omega _{\xi ,f}\otimes \id \otimes \omega
_{\eta ,g})(V_{23}^*V_{12})=\rho (\omega _{g,\eta })^*L(\omega _{\xi
,f})\,$.

\cqfd

\prop{Soit $\;f\in {\cal H}\,,\; \| f \| = 1\;$ tel que $\;V(f\otimes
f)=f\otimes f\,$. Notons $\;\psi \;$ l'\'etat $\;\psi =\omega _{f,f}\,$.
\hfill\break
 a) L'op\'erateur $\;L(\psi )\;$ est le projecteur sur $\;H^f\,$;
l'op\'erateur
$\;\rho (\psi )\;$ est le projecteur sur $\;H_f\,$. Ces deux projecteurs
commutent.\hfill\break
  b) Pour tout $\; \xi \in H_f\otimes H^f\,,\; V\xi =\xi \,$. \hfill\break
	c) On a $\;L(\psi )e=\theta _{f,f}e\,,\; \rho (\psi )\hat e=\theta
_{f,f}\hat e\,,\;\langle e,\hat e\rangle=
\langle e,f\rangle\langle f,\hat e\rangle\;$ et $\;Jf =\widehat
Jf=Uf=f\,$. \hfill\break
 d) On a $\;H^f\cap H_f=\C f\;$ et $\;\dim(H^f)\dim(H_f)=\dim({\cal
H})\,$. De plus $\; |\langle e,f\rangle |^2\dim(H_f)=|\langle \hat
e,f\rangle |^2\dim(H^f)=1\,$.}

\pf a) Par le lemme \sn 2.a), $\;L(\psi)\;$ et $\;\rho (\psi )\;$ sont des
idempotents. Comme $\;\| L(\psi ) \|\le 1\;$ et $\;\| \rho(\psi ) \|\le
1\,,$ ce sont des projecteurs. Les \'egalit\'es $\;L(\psi ){\cal H}=H^f\;$
et
$\;\rho(\psi ){\cal H}=H_f\;$ r\'esultent du lemme \sn1.a). Par le lemme
\sn 2.b)\ref, les projecteurs $\;L(\psi )\;$ et $\;\rho (\psi )=\rho(\psi
)^*\;$ commutent.

\noindent b) Pour $\;\xi \in H^f\,,\; \zeta \in H_f\,, \zeta \otimes f\otimes
\xi \;$ est invariant par $\;V_{12}\;$ et par $\;V_{23}\;$ donc par
$\;V_{13}\,$.

\noindent c) On a $\;\langle e,L(\psi )e\rangle =\varphi (L(\psi ))=\psi (\rho
(\varphi ))\,$. Or $\;\rho (\varphi )\;$ \'etant le projecteur sur l'espace
des vecteurs fixes, par hypoth\`ese r\'eduit \`a $\;\C e\,$, on en d\'eduit
$\;\langle e,L(\psi )e\rangle =|\langle f,e\rangle |^2=\langle e,\theta
_{f,f}e\rangle \,$. Comme $\;L(\psi )-\theta _{f,f}\;$ est positif, on a
$\;L(\psi )e=\theta _{f,f}e\,$. L'\'egalit\'e $\;\rho (\psi )\hat e=\theta
_{f,f}\hat e\;$ s'en d\'eduit en rempla\c cant $\;V\;$ par $\;\Sigma
V^*\Sigma \,$. De plus, comme l'image du projecteur $\;\rho
(\varphi)=\theta _{ e,e}\;$ est contenue dans $\;H_f\,$, on a $\;\rho
(\varphi)=\rho (\varphi )\rho (\psi )\,$, donc $\;\langle e,\hat e\rangle
e=\rho (\varphi )\hat e=\rho (\varphi )\rho (\psi )\hat e=\rho (\varphi
)(\langle f,\hat e\rangle f)=\langle e,f\rangle \langle f,\hat e\rangle
e\,$.

On a $\;JL(\psi )e=L(\psi)^*e=L(\psi )e\,$, $\;\widehat JL(\psi )e=
L(\psi^*)e=L(\psi )e\,$. Comme $\;e\;$ est s\'eparateur pour $\;S\,$,
$\;L(\psi )e\ne 0\,$; donc $\;Jf=\widehat Jf=f\,$, vu que $\;f\;$ est
proportionnel \`a $\;L(\psi )e\,$. Donc $\;Uf=\widehat JJf=f\,$

\noindent d) On a $\;\psi (\rho (\tau ))=\tau (L(\psi ))\,$. Mais $\;\rho
(\tau )=\theta _{e,e}\,$; donc $\;\psi (\rho (\tau ))= |\langle f,e\rangle
|^2\,$. D'autre part, $\;n\tau (L(\psi ))\;$ est la dimension de $\;H^f\;$
(o\`u
$\;n=\dim \H\,$). Donc $\;\dim H^f=n|\langle e,f\rangle|^2= n|\langle \hat
e,e\rangle|^2|\langle \hat e,f\rangle|^{-2}\;$ par c); donc
$\;|\langle \hat e,f\rangle |^2\dim(H^f)=1\,$. De m\^eme, $\; |\langle
e,f\rangle |^2 \dim(H_f)=1\,$. On en d\'eduit que $\;n^{-1}\dim H^f\dim
H_f=|\langle\hat e,e\rangle |^2\dim H^f\dim H_f=|\langle \hat e,f\rangle
|^2\dim(H^f) |\langle e,f\rangle |^2 \dim(H_f)=1\,$. Or par le corollaire.
1.8\ref, $\;\dim H_f\cap H^f=n^{-1}\dim H_f\dim H^f\,$.

\cqfd

\df{On appelle {\it pr\'e-sous-groupe}\/ de $\;V\;$ tout vecteur
$\;f\in {\cal H}\,,$ tel que $\;\|f\|=1\,,\; \langle f,e\rangle > 0
\,,\; V(f\otimes f)=f\otimes f\,$. Si $\;f\;$ est un pr\'e-sous-groupe,
les projecteurs $\;L(\omega _{f,f})\;$ et $\;\rho (\omega _{f,f})\;$ se
notent simplement $\;L_f\;$ et $\;\rho _f\,$.}

Par hypoth\`ese $\;e\;$ et $\;\hat e\;$ sont des pr\'e-sous-groupes.

\medskip\noindent{\it Exemple. ---}\/ Soit $\;G\;$ un groupe fini.
Notons $\;V\in \L(\ell^2(G)\otimes \ell^2(G)\;$ son unitaire
multiplicatif, donn\'e par la formule $\;(V\xi )(s,t)=\xi (st,t)\,$, pour
tout $\;\xi \in \ell^2(G)\otimes \ell^2(G)=\ell^2(G\times G)\,,\;s,t\in
G\,$. Soit $\;\Gamma \;$ un sous-groupe de $\;G\,$; alors $\;|\Gamma
|^{-1/2}\chi _\Gamma \;$ est un pr\'e-sous-groupe de $\;V\,$, o\`u
$\;\chi _\Gamma \;$ d\'esigne la fonction caract\'eristique de
$\;\Gamma \,$. Inversement, soit $\;f\in\ell^2(G)\;$ un
pr\'e-sous-groupe; alors, $\;\|f\|=1\,$, $\;f(1)=\langle \hat e,f\rangle
\in \R_+\;$ et, pour tout $\;s,t\in G\,$, on a $\;f(s)f(t)=f(st)f(t)\,$;
alors, l'ensemble $\;\Gamma =\{\,s\in G\,,\;f(s)\ne 0\,\}\;$ est stable
par le produit de $\;G\,$; comme $\;G\;$ est fini (et $\;f\ne 0\,$),
$\;\Gamma \;$ est un sous-groupe de $\;G\,$; de plus, pour tout
$\;s,t\in \Gamma \,$, on a  $\;f(st)=f(s)\,$; donc $\;f=|\Gamma
|^{-1/2}\chi _\Gamma \,$. On voit alors que $\;H^f\;$ est l'ensemble
des fonctions support\'ees par $\;\Gamma \;$ et $\;H_f\;$ est
l'ensemble des fonctions invariantes \`a droite par $\;\Gamma \,$. La
relation $\;\dim(H^f)\dim(H_f)=\dim(\H)\;$ (prop. \sn3.d) est le
th\'eor\`eme de Lagrange.

\prop{Soient $\;f,g\;$ des pr\'e-sous-groupes.\hfill\break
 a) On a $\;L(\omega _{f,g})^2=\langle f,g\rangle L(\omega _{f,g})\;$ et
$\;\rho (\omega _{f,g})^2= \langle f,g\rangle \rho (\omega _{f,g})\,$. On
a $\;L(\omega _{f,g})=L(\omega _{f,g})^*\;$ et $\;\rho (\omega
_{f,g})=\rho (\omega _{f,g})^*\,$.\hfill\break
 b) Les op\'erateurs $\;L(\omega _{f,g})\;$ et $\;\rho _g\;$ commutent.
Les op\'erateurs $\;\rho (\omega _{f,g})\;$ et $\;L _f\;$
commutent. \hfill\break
 \indent On note $\;H^{f,g}\;$ l'image de $\;L(\omega _{f,g})\;$ et
$\;H_{g,f}\;$ celle de $\;\rho(\omega _{g,f})\,$.\hfill\break
 c) On a $\;\dim(H^{f,g}\cap H_g)\dim (H^g)=\dim( H^{f,g})\;$ et
$\;\dim(H_{g,f}\cap H^g)\dim (H_g)=\dim( H_{g,f})\;$.  \hfill\break
 d) On a $\;\dim(H^{f,g}\cap H_g)\langle f,g\rangle\langle
g,e\rangle=\langle f,e\rangle\;$ et $\;\dim(H_{g,f}\cap H^g)\langle
f,g\rangle\langle f,e\rangle=\langle g,e\rangle\,$; en particulier,
$\;\langle f,g\rangle > 0\,$.\hfill\break
 e) On a $\;\dim(H^{f,g}\cap H_g)\dim (H_{g,f}\cap H^g)=\langle f,g
\rangle ^{-2}\,$.\hfill\break
 f) On a $\;H^f\subset H^{f,g}\,,\;H^g\subset H^{f,g}\,,\; H_f\subset
H_{g,f}\;$ et $\;H_g\subset H_{g,f}\,$.}

\pf a) Les deux premi\`eres assertions r\'esultent du lemme \sn2.a); les
deux derni\`eres de la prop. \sn.3.c) et de 1.6\ref.

\noindent b) r\'esulte du lemme \sn2.b).

\noindent c) Il r\'esulte du corollaire 1.8\ref, et de b) que
$\;n\dim(H^{f,g}\cap H_g) =\dim( H^{f,g})\dim H_g\,$. Or $\;n=\dim
H_g\dim H^g\,$. La deuxi\`eme assertion se d\'emontre de mani\`ere
analogue.

\noindent d) De c) on d\'eduit que $\;\langle f,g\rangle \dim(H^{f,g}\cap
H_g)\dim (H^g)=\langle f,g\rangle \dim H^{f,g}=n\tau (L(\omega_{f,g}))
=n\omega_{f,g}(\rho(\tau))=n\langle f,e\rangle\langle g,e\rangle\,$,
donc
$\;\langle f,g\rangle \dim(H^{f,g}\cap H_g)=\dim H_g\langle f,e\rangle
\langle g,e\rangle\,$. On en d\'eduit la premi\`ere assertion vu que
$\;\langle g,e\rangle^2\dim H_g=1\,$. La deuxi\`eme assertion se
d\'emontre de mani\`ere analogue, vu que $\;\langle f,e\rangle \langle
f,\hat e\rangle= \langle g,e\rangle \langle g, \hat e\rangle \,$.

\noindent e)  r\'esulte imm\'ediatement de d).

\noindent f) On a $\;L(\omega_{f,g})L_fe=\langle f,e\rangle
L(\omega_{f,g})^*f=\langle f,e\rangle \langle f,g\rangle f=\langle
f,g\rangle L_fe\,$. De m\^eme, $\;L(\omega_{f,g})L_ge=\langle g,e\rangle
L(\omega_{f,g})g=\langle g,e\rangle \langle f,g\rangle g=\langle
f,g\rangle L_ge\,$. On en d\'eduit que $\;L(\omega_{f,g})L_f =\langle
f,g\rangle L_f\;$ et $\;L(\omega_{f,g})L_g=\langle f,g\rangle L_g\,$, vu
que
$\;e\;$ est s\'eparateur pour $\;S\,$, d'o\`u les premi\`eres assertions,
vu que $\;\langle f,g\rangle \ne0\;$ (par d). Les autres assertions s'en
d\'eduisent en rempla\c cant $\;V\;$ par $\;\Sigma V^*\Sigma\,$.

\cqfd

\cor{Il y a un nombre fini de pr\'e-sous-groupes pour $\;V\,$.}

\pf Il r\'esulte imm\'ediatement de la proposition pr\'ec\'edente que, si
$\;f\;$ et  $\;g\;$ sont des pr\'e-sous-groupes distincts, $\;\langle
f,g\rangle ^{-2}\;$ est entier, donc $\;\langle f,g\rangle \le 2^{-1/2}\,$,
donc $\;\|f-g\|^2\ge 2-\sqrt 2\,$, d'o\`u le r\'esultat.

\cqfd

\prop{a) Soient $\;f\;$ et $\;g\;$ des pr\'e-sous-groupes. Les conditions
suivantes sont \'equivalentes.\hfill\break
 (i) $\;V(f\otimes g)=f\otimes g\,$. \quad (ii) $\;H^{g}\subset H^f\,$.
\quad (iii) $\;H_f\subset H_{g}\,$. \quad (iv) $\;H^{f,g}=H^f\,$; \
 (v) $\;H^{g,f}=H^f\,$; \quad
 (vi) $\;H_{f,g}=H_g\,$; \quad
 (vii) $\;H_{g,f}=H_g\,$; \hfill\break
 \indent Si $\;f\;$ et $\;g\;$ v\'erifient ces conditions, nous \'ecrirons
$\;g\prec f\,$.\hfill\break
 b) La relation $\;\prec\;$ est une relation d'ordre sur l'ensemble des
pr\'e-sous-groupes de $\;V\,$, pour laquelle $\;e\;$ est le plus grand
\'el\'ement et $\;\hat e\;$ est le plus petit \'el\'ement.}

\pf a) (ii)$\;\Ra\;$(i) r\'esulte de ce que $\;g\in H^{g}\,$. Si (i) est
v\'erifi\'ee, pour tout $\;\xi \in H^{g}\,$, on a
$\;V(f\otimes \xi)= f\otimes \xi\;$ par la prop. \sn3.b), donc (i)$\Ra$(ii).
(i)$\iff $(iii) en r\'esulte rempla\c cant $\;V\;$ par $\;\Sigma
V^*\Sigma
\,$. Par la prop. \sn 5.f), (iv)$\Ra$(ii), (v)$\Ra$(ii), (vi)$\Ra$(iii) et
(vii)$\Ra$(iii). Si (i) est satisfaite, par le lemme \sn2.a), $\;L_fL(\omega
_{f,g})=L(\omega _{f,g})\,,\;L_fL(\omega _{g,f})=L(\omega_{g,f})\,,\; \rho
(\omega _{f,g})\rho _g=\rho (\omega _{f,g})\,,\;\rho (\omega _{g,f})\rho
_g=\rho (\omega_{g,f})\,,$ donc (i) implique les relations (iv) \`a (vii).

\noindent b) Il est clair, par la condition (ii) que $\;\prec \;$ est une
relation de pr\'eordre. Si $\;f\prec g\;$ et $\;g\prec f\,$, on a $\;g\in
H^f\cap H_f=\C f\,$. Comme $\;\|f\|=\|g\|=1\,,\;\langle e,f\rangle >0\;$
et
$\;\langle e,g\rangle >0\,,\; f=g\,$. Il est clair enfin que $\;e\;$ est le
plus grand \'el\'ement et $\;\hat e\;$ le plus petit.

\cqfd

\prop{Soient $\;f\;$ et $\;g\;$ deux pr\'e-sous-groupes. Si le nombre
entier $\;\langle f,g\rangle^{-2}\;$ est premier alors $\;f\;$ et $\;g\;$
sont comparables.}

\pf Par la prop. \sn5.e) un des nombres $\;\dim(H^{f,g}\cap H_g)\;$ ou
$\;\dim (H_{g,f}\cap H^g)\;$ est \'egal \`a 1; le r\'esultat d\'ecoule
alors des prop. \sn5.c) et \sn7.a) (conditions (v) et (vii)).

\cqfd

 \prop{Soient $\;f\,,\,g\;$ des pr\'e-sous-groupes tels que $\;g\prec
f\,$.\hfill\break
 a) $\;\rho _gL_f=L_f\rho _g\;$ est le projecteur sur
$\;H_{g}\cap H^f\,$. \hfill\break
 b) On a $\;L_gf=\langle g,f\rangle g\,,\; \langle g,e\rangle =\langle
g,f\rangle \langle f,e\rangle \,$. \hfill\break
 c) On a $\;\dim(H^f)\dim(H_{g})=n\dim(H_{g}\cap H^f)\;$ et
$\;\dim(H_{g}\cap H^f)\langle f,g\rangle ^2=1\,$.}

\pf a) r\'esulte du lemme \sn 2.b)\ref.

\noindent b) On a $\;L_fe=\theta _{f,f}e=\langle f,e\rangle f\,$. De m\^eme
$\;L_ge=\langle g,e\rangle g\,$. Comme $\;L_g=L_gL_f\,$, on en d\'eduit
que $\;\langle f,e\rangle L_gf=\langle g,e\rangle g\,$. Comme $\;L_gf\;$
est colin\'eaire \`a $\;g\,$, $\;L_gf=\langle g,f\rangle g\,$. On a alors
$\;\langle g,f\rangle \langle f,e\rangle g=\langle f,e\rangle
L_gf=L_ge=\langle g,e\rangle g\,$.

\noindent c) La premi\`ere assertion se d\'eduit imm\'ediatement du corollaire
1.8\ref. La deuxi\`eme de la prop. \sn5.e) et de \sn3.d), vu que
$\;H^{f,g}=H^f\;$ et $\;H_{g,f}=H_g\,$.

 \cqfd

\lem{Soit $\;\psi \;$ un \'etat sur $\;{\cal L}({\cal H})\;$ tel que
$\;L(\psi )\;$ soit un projecteur. Alors $\;V(L(\psi )e\otimes L(\psi )e)=
L(\psi )e\otimes L(\psi )e\,$.}

\pf La restriction de $\;\psi \;$ \`a $\;\widehat S\;$ est de la forme
$\;x\widehat \varphi \,$, o\`u $\;x\;$ est un \'el\'ement positif de
$\;\widehat S\,$. Posons $\;\xi =x^{1/2}\hat e\,$. Les restrictions de
$\;\psi \;$ et $\;\omega _{\xi ,\xi }\;$ \`a $\;\widehat S\;$ co\"incident,
donc $\;L(\psi )=L(\omega _{\xi ,\xi })\,$. On peut donc supposer que
$\;\psi =\omega _{\xi ,\xi }\,$. Par le lemme \sn1.a),
$\;L(\psi ){\cal H}=H^\xi \,$.

Posons $\;\eta = L(\psi )e\,$; par le lemme \sn1.b), $\;V(\xi \otimes
e)\in {\cal H}\otimes H_\eta \,$, d'o\`u $\;\eta =L(\omega _{\xi ,\xi
})e\in H_\eta\,$.

\cqfd

\prop{a) Pour un projecteur orthogonal non nul $\;p\in S\;$ les conditions
suivantes sont \'equivalentes\quad  (i) $\;p\otimes p\le \delta
(p)\,$.\quad (ii) $\;\delta (p)(1\otimes p)=p\otimes p\,$.\quad (iii)
$\;\delta (p)(p\otimes 1)=p\otimes p\,$.\quad
 (iv) $\;\rho(\omega _{e,pe})\ge 0\,$.\quad
 (v) $\;\varphi (p)^{-1}\rho(\omega _{e,pe})\;$ est un projecteur
orthogonal.
\quad (vi) $\;f=\|pe\|^{-1}pe\;$ est  un pr\'e-sous-groupe et
$\;p=L_f\,$.\hfill\break  b) Pour un projecteur orthogonal non nul $\;p\in
\widehat S\;$ les  conditions suivantes sont \'equivalentes. \quad
 (i) $\;p\otimes p\le \widehat\delta (p)\,$.\quad
 (ii) $\;\widehat\delta (p)(1\otimes p)=p\otimes p\,$.\quad
 (iii) $\;\widehat\delta (p)(p\otimes 1)=p\otimes p\,$.\quad
 (iv) $\;L(\omega _{\hat e,p\hat e})\ge 0\,$.\quad
 (v) $\;\widehat\varphi (p)^{-1}L(\omega _{\hat e,p\hat e})\;$ est un
projecteur orthogonal.\quad (vi)$\;f=\|p\hat e\|^{-1}p\hat e\;$ est  un
pr\'e-sous-groupe et $\;p=\rho_f\,$.}

\pf Montrons a). Rempla\c cant $\;V\;$ par $\;\Sigma V^*\Sigma\,,$ on
en d\'eduit b).

\noindent (i) $\Ra$ (vi). On a $\;(p\otimes p)V(pe\otimes pe)=(p\otimes
p)\delta (p)(1\otimes p)(e\otimes e)=pe\otimes pe\,$, d'o\`u il ressort,
comme $\;\|V(pe\otimes pe)\|=\|pe\otimes pe\|\;$ que $\;V(pe\otimes
pe)=pe\otimes pe\,$. Par ailleurs $\;pe\ne 0\,$, puisque $\;e\;$ est
s\'eparateur pour $\;S\,$. Il s'ensuit que $\;f=1/\|pe\|pe\;$ est un
pr\'e-sous-groupe. On a $\;L_f(e)=\langle f,e\rangle f=pe\,$, donc
$\;L_f=p\,$.

\noindent (vi)$\Ra$(v). Si $\;p=L_f\,$, $\;\rho (\omega _{e,pe})=\langle
e,f\rangle \rho (\omega _{e,f})=\langle e,f\rangle^2\rho _f\;$ par la
prop. \sn 7.a), condition (vi).

\noindent (v)$\Ra$(iv) est clair.

\noindent (iv)$\Ra$(ii). Si (iv) est v\'erifi\'ee, il existe une forme
positive
$\;\psi \;$ telle que $\;p=L(\psi )\;$ (prop. 1.4). Comme $\;p\ne 0\,,\;\psi
(1)\ne 0\,$; or $\;\varepsilon :L(\omega )\mt \omega (1)\;$ est un
caract\`ere sur $\;S\,$; donc $\;\psi (1)=\varepsilon (p)=1\,$. Par le
lemme \sn 10, $\;(\delta (p)(1\otimes p)(e\otimes e)=V(pe\otimes
pe)=pe\otimes pe\,$, d'o\`u (ii), vu que $\;e\otimes e\;$ est
s\'eparateur pour $\;S\otimes S\,$.

\noindent Il est clair que (ii) implique (i).

En changeant $\;V\;$ en $\;(U\otimes 1)V^*(U\otimes 1)\,$, on change
$\;\delta \;$ en $\;\sigma \circ\delta \,$. De (i)$\iff $(ii), on d\'eduit
que (i)$\iff$(iii).

\cqfd

\cor{Soient $\;p\in S\;$ et $\;q\in \widehat S\;$ des projecteurs tels
que $\;p\hat e\ne 0\,,\;qe\ne 0\,,\;pq=qp\;$ et $\;\tau (pq)=1/n\,$.
Alors $\;\|pe\|^{-1}pe=\|q\hat e\|^{-1}q\hat e\,$; c'est un
pr\'e-sous-groupe $\;f\;$ et l'on a $\;p=L_f\;$ et $\;q=\rho _f\,$.}

\pf Comme $\;\theta _{\hat e,\hat e}\;$ est central dans $\;S\,$, on a
$\;p\hat e=\hat e\,$; de m\^eme,  $\;qe=e\,$. On a alors $\;\langle
pe,q\hat e\rangle =\langle e,pq\hat e\rangle =\langle qe,p\hat
e\rangle =n^{-1/2}\,$. Or $\;\|pe\|^2=\tau (p)\;$ et $\;\|q\hat e\|^2=\tau
(q)\,$; donc $\;\|pe\|\|q\hat e\|=n^{-1/2}\,$. Donc $\;pe\;$ et  $\;q\hat
e\;$ sont (positivement) proportionnels. Alors $\;p=n^{1/2}L(\omega
_{\hat e,pe})\;$ est proportionnel \`a $\;L(\omega _{\hat e,q\hat e})\,$,
d'o\`u l'on d\'eduit que $\;q\;$ est de la forme  $\;\rho _f\;$ (prop. \sn
11.b), condition (iv)). Comme $\;L_fe\;$ est proportionnel \`a $\;f\,$, on
en d\'eduit que $\;p=L_f\,$.

\cqfd

\prop{a) Soit $\;\psi \;$ un \'etat sur $\;{\cal L}({\cal H})\,$.  Il existe
des pr\'e-sous-groupes $\;f,g\;$ tels que $\;\{\,\xi \in{\cal H}\,,\;L(\psi
)\xi =\xi \,\}=H^f\;$ et $\;\{\,\xi \in{\cal H}\,,\;\rho (\psi )\xi =\xi
\,\}=H_g\,$.\hfill\break
 b) Soit $\;K\;$ un sous-espace non nul de $\;{\cal H}\,$. Il existe des
pr\'e-sous-groupes $\;f,g\;$ tels que $\;\{\,\xi \in {\cal H}\,,\; V(\zeta
\otimes \xi )=\zeta \otimes \xi \,,\; \forall \zeta \in K\,\}=H^f\;$ et
$\;\{\,\xi \in {\cal H}\,,\; V(\xi \otimes \zeta )=\xi \otimes \zeta \,,\;
\forall \zeta \in K\,\}=H_{g}\,$.}

\pf a) Notons $\;C=\{\,x\in S\,,\;\rho (\omega _{e,xe})\ge 0\,\}\,$. Par
la prop. 1.4.a\ref), $\;(1+L(\psi ))^k2^{-k}\in C\,$, donc sa limite, qui
est le projecteur orthogonal $\;p\;$ sur $\;\{\,\xi \in{\cal H}\,,\;L(\psi
)\xi =\xi \,\}\;$ aussi. L'assertion sur $\;f\;$ r\'esulte alors de la prop.
\sn11.a) (condition (iv)). L'assertion sur $\;g\;$ s'en d\'eduit en
rempla\c cant $\;V\;$ par $\;\Sigma V^*\Sigma \,$.

\noindent b) Il suffit d'appliquer a) \`a l'\'etat $\;\psi =k^{-1}\sum_i
\omega _{\zeta_i,\zeta _i}\;$ o\`u $\;\{\,\zeta _i\,,\; i=1,...,k\,\}\;$ est
une base orthonormale de $\;K\,$.

\cqfd

\cor{L'ensemble des pr\'e-sous-groupes ordonn\'e par $\;\prec \;$ est un
treillis. Plus pr\'ecis\'ement, \'etant donn\'e des pr\'e-sous-groupes
$\;f\;$ et $\;f'\,,$ il existe des pr\'e-sous-groupes $\;g\;$ et $\;g'\;$
tels que $\;H^f\cap H^{f'}=H^g\;$ et $\;H_f\cap H_{f'}=H_{g'}\,$.}

\pf Appliquons la proposition \sn 13.b) \`a $\;K=\C f+\C f'\,$. Soient
$\;g,g'\;$ des pr\'e-sous-groupes tels que $\;H^f\cap H^{f'}=H^g\;$ et
$\;H_f\cap H_{f'}=H_{g'}\,$; il est clair qu'alors $\;g=\inf(f,f')\;$ et
$\;g'=\sup(f,f')\,$.

\cqfd

Rappelons qu'une repr\'esentation de $\;V\;$ dans un espace hilbertien
$\;K\;$ est donn\'ee par un unitaire $\;X\in {\cal L}(K\otimes \H)\;$ tel
que $\;X_{12}X_{13}V_{23}=V_{23}X_{12}\in {\cal L}(K\otimes
\H\otimes \H)\,$. Comme toute repr\'esentation est un sous-multiple
de la r\'eguli\`ere, on a \'evidemment:

\cor{Soient $\;X\in {\cal L}(K\otimes \H)\;$ une repr\'esentation de
$\;V\;$ et $\;\xi \in K\,$. Il existe un pr\'e-sous-groupe $\;f\;$ tel que
$\;\{\eta\in\H\, ,\;X(\xi \otimes \eta)=(\xi\otimes \eta)\}=H^f\,$.}

\cqfd

\prop{Soit $\;\psi \;$ un \'etat de $\;{\cal L}({\cal H})\,$.\hfill\break
 a) Les ensembles des valeurs propres de module 1 de $\;L(\psi )\;$ et de
$\;\rho (\psi )\;$ sont des sous-groupes de $\;U(1)\,$; en particulier les
valeurs propres de module 1 de $\;L(\psi )\;$ et de
$\;\rho (\psi )\;$ sont des racines de l'unit\'e.\hfill\break
 b) Il existe des pr\'e-sous-groupes $\;f\;$ et $\;g\;$ tels que l'espace
engendr\'e par les vecteurs propres de $\;L(\psi )\;$ (resp. $\;\rho (\psi
)\,$) associ\'es aux valeurs propres de module 1 soit $\;H^f\;$ (resp.
$\;H_g\,$).}

\pf Il suffit d'\'etablir les \'enonc\'es relatifs \`a $\;L(\psi )\,$; les
\'enonc\'es relatifs \`a $\;\rho (\psi )\;$ s'en d\'eduisent en rempla\c
cant $\;V\;$ par $\;\Sigma V^*\Sigma \,$.

\noindent a) La restriction de $\;\psi \;$ \`a $\;\widehat S\;$
co\"incide avec celle d'un \'etat vectoriel $\;\omega _{\zeta ,
\zeta}\,$. Un nombre complexe $\;\lambda \;$ de module 1 est valeur
propre de $\;L(\psi )\;$ si et seulement s'il existe $\;\xi \in {\cal H}\;$
non nul tel que $\;V(\zeta \otimes \xi )=\lambda (\zeta \otimes \xi
)\,$. Soit $\;G\;$ l'ensemble des valeurs propres de module 1 de
$\;L(\psi )\,$. Si $\;\lambda ,\lambda '\in G\,$, il existe $\;\xi \;$ et
$\;\eta\;$ non-nuls avec $\;V(\zeta \otimes \xi )=\lambda '(\zeta
\otimes \xi )\;$ et $\;V(\zeta\otimes \eta)=\lambda (\zeta \otimes
\eta)\,$; on a alors $\;V_{12}V_{13}(\zeta \otimes \xi \otimes
\eta)=\lambda \lambda '(\zeta \otimes \xi \otimes \eta)\;$ donc
$\;\lambda\lambda '\;$ est valeur propre de $\;(\psi \otimes id\otimes
id)(V_{12}V_{13})\;$ donc de $\;L(\psi )\otimes 1\;$ qui lui est
conjugu\'e. De plus $\;1\;$ est valeur propre de vecteur propre
associ\'e $\;\hat e\,$. Comme $\;G\;$ est un sous-semi-groupe compact
(fini) de $\;U(1)\,$, c'est un sous-groupe.

\noindent b) Soit $\;\psi '\;$ l'\'etat de $\;{\cal L}({\cal H})\;$ tel que
$\;L(\psi ')=L(\psi )^k\;$ o\`u $\;k\;$ est l'ordre du groupe cyclique
$\;G\,$. L'espace engendr\'e par les vecteurs propres de $\;L(\psi )\;$
associ\'es aux valeurs propres de module 1 est l'espace propre de
$\;L(\psi ')\;$ pour la valeur propre 1; il r\'esulte de la proposition \sn
13.a) qu'il est de la forme $\;H^f\,$.

\cqfd

\sec{Pr\'e-sous-groupes et sous-alg\`ebres co\"id\'eales.}

Dans ce paragraphe, nous  allons \'etablir une correspondance naturelle
entre les pr\'e-sous-groupes de $\;V\;$ et les sous-alg\`ebres co\"id\'eales
(\cf \eg [\ref\Mon]) de la $\,C^*$-alg\`ebre de Hopf associ\'ee $\;S\,$. Nous
donnons ensuite dans notre cadre une d\'emonstration tr\`es simple de
r\'esultats de [\ref\Ma], ainsi que d'autres r\'esultats analogues.

Commen\c cons  par rappeler la d\'efinition des sous-alg\`ebres
co\"id\'eales dans une alg\`ebre de Hopf.

\df{ ([\ref\Mon]) Soit $\;(A,\delta)\;$ une alg\`ebre de Hopf unif\`ere. Une
sous-alg\`ebre  $\;B\;$ de $\;A\;$ contenant l'unit\'e est dite {\it
co\"id\'eale \`a droite}\/ (resp. {\it \`a gauche}) si $\;\delta(B)
\subset B\otimes A\;$ (resp. $\;\delta(B) \subset A\otimes B\,)$.}

Revenons \`a pr\'esent aux alg\`ebres de Hopf associ\'ees \`a notre unitaire
multiplicatif $\;V\;$ de multiplicit\'e $\;1\,$. Si $\;f\;$ est un
pr\'e-sous-groupe on pose $\;R_f=UL_fU\;$ et $\;\lambda _f=U\rho _fU\,$.
Posons aussi $\;G_f=\{\,x\in
S\,,\;xe\in H_f\,\}\,,
\;D_f=\{\,x\in S\,,\;xe\in UH_f\,\}\,, \;\widehat D_f= \{\,x\in \widehat
S\,,\;x\hat e\in H^f\,\}\;$ et $\;\widehat G_f=\{\,x\in \widehat
S\,,\;x\hat e\in UH^f\,\}\,$.

\prop{Soit $\;f\;$ un pr\'e-sous-groupe.\hfill\break
 a) On a $\;G_f=\{\,L(\omega _{\xi ,f})^*\,,\;\xi \in \H\,\}=\{\,x\in S\,,\;
[x,\rho  _f]=0\,\}=\{\,x\in S\, ,\; \delta (x)(1\otimes L_f)=x\otimes
L_f\,\}\,.$ C'est une sous-alg\`ebre co\"id\'eale \`a gauche de $\;S\,$, stable
par l'involution.\hfill\break
 b) On a $\;D_f=\{\,L(\omega _{f,\xi })\,,\;\xi \in {\cal H}\,\}=\{\,x\in
S\,,\; [x,\lambda _f]=0\,\}=\{\,x\in S\, ,\; \delta (x)(L_f\otimes
1)=L_f\otimes x\,\}\,.$ C'est une sous-alg\`ebre co\"id\'eale \`a droite de
$\;S\,$, stable par l'involution.\hfill\break
 c) On a $\;\widehat D_f=\{\,\rho (\omega _{f,\xi })\,, \;\xi \in {\cal
H}\,\}=\{\,x\in \widehat S\,,\; [x,L _f]=0\,\}=\{\,x\in \widehat S\, ,\;
\hat\delta (x)(\rho _f\otimes 1)= \rho _f\otimes x\,\}\,.$ C'est une
sous-alg\`ebre co\"id\'eale \`a droite de $\;\widehat S\,$, stable par
l'involution. \hfill\break
 d) On a $\;\widehat G_f=\{\,\rho (\omega _{\xi ,f })^*\,, \;\xi \in {\cal
H}\,\}=\{\,x\in \widehat S\,,\; [x,R _f]=0\,\}=\{\,x\in \widehat S\,, \;
\hat\delta (x)(1\otimes \rho _f)=x \otimes \rho _f\,\}\,.$ C'est une
sous-alg\`ebre co\"id\'eale \`a gauche de $\;\widehat S\,$, stable par
l'involution.}

\pf a) Par le lemme 3.2.b\ref), on a $\;\{\,L(\omega _{\xi ,f})^*\,,\;\xi
\in {\cal H}\,\}\subset \{\,x\in S\,,\; [x,\rho  _f]=0\,\}\,$. Si $\;x\rho
_f=\rho_f x\,$, on a $\;xe=x\rho_fe\in H_f\,$, donc $\;x\in G_f\,$. Soit
$\;x\in G_f\,$. Alors, par 1.3\ref, $\;\kappa(x)=n^{1/2}L(\omega _{\hat
e,Uxe})\,$, donc $\;x=n^{1/2}L(\omega _{Uxe, \hat e})^*\,$; or $\;L(\omega
_{Uxe,\hat e})^*=(\omega _{\hat e,Uxe}\otimes
\id )(V^*)=(\omega _{\hat e,U\rho_f xe}\otimes \id )(V^*)=(\omega
_{\hat e,Uxe}\otimes \id )(V^*(\lambda_f \otimes 1))=(\omega _{\hat
e,Uxe}\otimes \id )((\lambda_f \otimes 1)V^*) =L(\omega _{Uxe, \lambda
_f\hat e})^*\,$. Or $\;\lambda_f\hat e\;$ est proportionnel \`a $\;f\,$,
donc $\;x\in \{\,L(\omega _{\xi ,f})^*\,,\;\xi \in {\cal H}\,\}\,$.  Enfin,
pour $\;x\in S\,,\;\delta (x)(1\otimes L_f) =x\otimes L_f\iff \delta
(x)(1\otimes L_f)(e\otimes e)=xe\otimes L_f e\iff V(xe\otimes f)=xe
\otimes f\iff xe\in H_f\,$.

Il est clair que $\;G_f=\{\,x\in S\,,\; [x,\rho  _f]=0\,\}\;$ est une
sous-alg\`ebre involutive unitale de $\;S\,$. Si $\;x\in G_f\,,$ on a
$\;\delta (x)(e\otimes e)=V(x\otimes 1)(e\otimes e)\in V(H_f\otimes
H_f)\subset \H\otimes H_f)\;$ par le lemme 3.1\ref; donc $\;\delta
(x)\in S\otimes G_f\,$.

\medskip\noindent b), c) et d) se d\'eduisent de a) en rempla\c cant
successivement
$\;V\;$ par $\;(U\otimes 1)V^*(U\otimes 1)\,$, $\;\Sigma V^*\Sigma\;$ et
$\;\widetilde V\,$.

\cqfd

Notons que $\;G_f=\kappa (D_f)=\widehat J D_f\widehat J\;$ et $\;\widehat
G_f=\widehat\kappa (\widehat D_f)=J\widehat D_f J\,$.

\prop{a) L'application $\;f\mt D_f\;$ (resp. $\;f\mt G_f\,$) est une
bijection entre pr\'e-sous-groupes et sous-alg\`ebres co\"id\'eales \`a
droite (resp \`a gauche) de l'alg\`ebre de Hopf $\;S\,$. L'application $\;f\mt
\widehat D_f\;$ (resp. $\;f\mt\widehat G_f\,$) est une bijection entre
pr\'e-sous-groupes et sous-alg\`ebres co\"id\'eales \`a droite (resp \`a
gauche) de l'alg\`ebre de Hopf $\;\widehat S\,$.\hfill\break
 b) Toute sous-alg\`ebre co\"id\'eale (\`a gauche ou \`a droite) de $\;S\;$
est stable par l'involution de $\;S\,$.}

\pf a) Soit $\;B\;$ une sous-alg\`ebre co\"id\'eale \`a droite de $\;S\,$.
Posons $\;H=Be\,$. Par hypoth\`ese, pour $\;x\in B\;$ et $\;\xi\in \H\,$, on a
$\;V(xe\otimes \xi) =\delta (x)(e\otimes \xi)\in H\otimes \H\,$. Donc
$\;V(H\otimes \H)=H\otimes \H\,$. Notons $\;p\;$ le projecteur sur $\;H\,$.
Le projecteur $\;p\otimes 1\;$ commute \`a $\;V\,$, donc $\;p\in U\widehat
SU\,$. De plus, pour $\;x,y\in B\,$, on a $\;\widehat V^*(xe\otimes
ye)=\delta (y)(xe\otimes e)\in H\otimes \H\,$. Donc $\;H\otimes
H\subset \widehat V(H\otimes \H)\,$, donc $\;p\otimes p\le \widehat
V(p\otimes 1)\widehat V^*\,$. Il r\'esulte alors de la prop. 3.11 \ref
(appliqu\'ee \`a $\;\widehat V\,$), qu'il existe un pr\'e-sous-groupe
$\;f\;$ tel que $\;p=\lambda _f\,$, donc $\;B=D_f\,$.

Si $\;B\;$ est une sous-alg\`ebre co\"id\'eale \`a gauche, $\;\kappa (B)\;$
est une sous-alg\`ebre co\"id\'eale \`a droite; l'assertion resp. s'en d\'eduit.
Les autres assertions s'en d\'eduisent en rempla\c cant $\;V\;$ par
$\;\widetilde V\,$.

\noindent b) se d\'eduit imm\'ediatement de a) et de la prop. \sn2.

\cqfd

Comme $\;e\;$ est un vecteur totalisateur et s\'eparateur pour
$\;S\,$, on a $\;\dim D_f=\dim UH_f=\dim H_f=\dim G_f\,$; comme
$\;\hat e\;$ est un vecteur totalisateur et s\'eparateur pour $\;\widehat
S\,$, on a $\;\dim \widehat D_f=\dim H^f=\dim UH^f=\dim \widehat
G_f\,$. On en d\'eduit alors \`a l'aide de la prop. 3.3.d\ref):

\cor{(Th\'eor\`eme de Lagrange, \cf [\ref\NZ, \Ma]) La dimension de toute
sous-alg\`ebre co\"id\'eale divise la dimension de $\;S\,$.}

\cqfd

\prop{Soient $\;f\;$ et $\;g\;$ des pr\'e-sous-groupes. L'espace vectoriel
engendr\'e par $\;\{xy\,,\;x\in D_f\,,\;y\;\in \widehat G_g\,\}\;$ est une
sous-alg\`ebre involutive de $\;{\cal L} ({\cal H})\,$.}

\pf On a $$\eqalign{\rho (\omega _{\eta ,g})^*L(\omega _{f,\xi})
&=(\omega _{f,\xi }\otimes \id\otimes \omega
_{g,\eta})(V_{23}^*V_{12})=(\omega _{f,\xi }\otimes \id\otimes \omega
_{g,\eta})(V_{12}V_{23}^*V_{13}^*)\cr&=\sum_i L(\omega
_{f,\xi_i})\rho (\omega _{\eta _i,g})^*\,,}$$ o\`u on a pos\'e $\;V^*(\xi
\otimes \eta)=\sum_i \xi_i\otimes \eta_i\,$, d'o\`u le r\'esultat.

\cqfd

On note $\;B_{f,g}\;$ la sous-alg\`ebre d\'efinie dans la proposition \sn 5.
Rempla\c cant $\;V\;$ par $\;\widehat V\;$ on en d\'eduit que l'espace
vectoriel engendr\'e par $\;\{xy\,,\;x\in U\widehat D_gU\,,\;y\in G_f
\,\}\;$ est une sous-alg\`ebre involutive de $\;{\cal L}({\cal H})\,$,
not\'ee $\;A_{g,f}\,$.

\medskip\noindent{\it Remarque. ---}\/ On peut donner la
g\'en\'eralisation suivante \`a la prop. \sn5: soit $\;A\;$ une
$\,C^*$-alg\`ebre munie d'une coaction $\;\delta _A:A\ra A\otimes S\;$ de
$\;S\,$. Alors, avec les notations de [\BSb\ref], le sous-espace vectoriel de
$\;A\sd \widehat S\;$ engendr\'e par $\;\{\,\widehat \theta (x)\pi (a)\,,\;
a\in A\,,\; x\in \widehat G_g\,\}\;$ est une sous-alg\`ebre involutive de
$\;A\sd \widehat S\,$. En effet, pour $\;a\in A\;$ et $\;\xi \in \H\,$, on a
$\;\pi _L(a)(1\otimes \rho (\omega _{\xi ,g})^*)=(\id\otimes \id\otimes
\omega _{g,\xi })(\delta _A(a) _{12}V_{23}^*)=(\id\otimes \id\otimes
\omega _{g,\xi })(V_{23}^*(\id\otimes \delta )\circ\delta _A(a))=
(\id\otimes \id\otimes \omega _{g,\xi })(V_{23}^*(\delta _A\otimes \id
)\circ\delta _A(a))=\dst\sum_i(\id\otimes \id\otimes \omega _{g,s_i\xi
})(V_{23}^*\delta _A(a_i)_{12})\;$ o\`u l'on a pos\'e $\;\delta
_A(a)=\dst\sum_i a_i\otimes s_i\,$.

Remarquons que, comme $\;D_f\;$ est une sous-alg\`ebre co\"id\'eale \`a
droite de $\;S\,$, le coproduit de $\;S\;$ d\'efinit par restriction une
coaction de $\;S\;$ dans $\;D_f\,$; la prop. 4.5 est donc un cas particulier du
r\'esultat ci-dessus.

\lem{On a: $\;\dim B_{f,g}= \dim H_f\dim H^g=\dim A_{g,f}\,$.}

\pf L'application $\;\mu :x\otimes y\mt xy\;$ est une application lin\'eaire
bijective de $\;S\otimes \widehat S\;$ sur $\;{\cal L}({\cal H})\,$, d'o\`u
la premi\`ere assertion, vu que $\;B_{f,g}=\mu(D_f\otimes \widehat
G_g)\,$.

La deuxi\`eme assertion en r\'esulte en rempla\c cant $\;V\;$ par
$\;\widehat V\,$.

\cqfd

\lem{Soit $\;f\;$ un pr\'e-sous-groupe. Le commutant de $\;G_f\;$ (resp.
$\;D_f\,,\; \widehat D_f\,,\; \widehat G_f\,$) est $\;UA_{f,\hat e}U\;$ (resp.
$\;UB_{\hat e,f}U\,,\; A_{e,f}\,,\;UB_{f,e}U\,$).}

\pf Par le lemme 3.2.b\ref), ces deux alg\`ebres commutent. Par ailleurs
$\;G_f=\{\,x\in S\,,\; [x,\rho_f]=0\,\}\supset (USU)'\cap \widehat
D_f'=(UA_{f,\hat e}U)'\,$. Donc $\;G_f=(UA_{f,\hat e}U)'\,$; le lemme
r\'esulte du th\'eor\`eme du bicommutant.

On en d\'eduit les assertions \og resp.\fg\  en rempla\c cant successivement
$\;V\;$ par $\;(U\otimes 1)V^*(U\otimes 1)\,$, $\;\Sigma V^*\Sigma\;$ et
$\;\widetilde V\,$.

\cqfd

On peut imm\'ediatement g\'en\'eraliser ce lemme:

\prop{Soient $\;f\;$ et $\;g\;$ des pr\'e-sous-groupes. Le commutant de
$\;B_{f,g}\;$ (resp. $\;A_{f,g}\,$) est $\;UB_{g,f}U\;$ (resp.
$\;UA_{g,f}U\,$).}

\pf On a $\;B_{f,g}'=D_f'\cap \widehat G_g'=U(B_{\hat e,f}\cap B_{g,e})U\,$.
L'application $\;\mu :x\otimes y\mt xy\;$ est une application lin\'eaire
bijective de $\;S\otimes \widehat S\;$ sur $\;{\cal L}({\cal H})\,$; donc
$\;B_{\hat e,f}\cap B_{g,e}=\mu((D_g\otimes \widehat S)\cap (S\otimes
\widehat G_f))=B_{g,f}\,$. L'assertion \og resp.\fg\ en r\'esulte
en rempla\c cant $\;V\;$ par $\;\widehat V\,$.

\cqfd

Le lemme \sn 7 nous permet de donner, dans notre cadre, une d\'emonstration
imm\'ediate de la g\'en\'eralisation par Masuoka du th\'eor\`eme de
Nichols-Zoeller (\cf [\ref\NZ, \Ma]).

\th {Soit $\;B\;$ une sous-alg\`ebre co\"id\'eale (\`a gauche ou \`a droite)
de $\;S\,$. Alors $\;S\;$ est un $\,B$-module (\`a gauche) libre.}

Remarquons que, gr\^ace \`a l'involution, on en d\'eduit imm\'ediatement
que $\;S\;$ est un $\,B$-module \`a droite libre.

\medskip\pf L'application $\;x\mt xe\;$ est un isomorphisme de $\,B$-modules
entre $\;S\;$ et $\;\H\,$. Le th\'eor\`eme se d\'eduit facilement du lemme
\sn 7 et du lemme \'el\'ementaire suivant.

\lem{Soient $\;H\;$ un espace hilbertien de dimension finie et $\;A\;$ une
sous-alg\`ebre involutive unitale de $\;{\cal L}(H)\,$. Les conditions
suivantes sont \'equivalentes:\hfill\break
 (i) Il existe $\;k\in \N-\{0\}\;$ tel que $\;H^k\;$ soit un $\,A$-module
libre.\hfill\break
 (ii) $\;\dim A\dim A'=(\dim H)^2\,$.\hfill\break
Dans ce cas, $\;H\;$ est un $\;A\;$ module libre si et seulement si
$\;\dim A\;$ divise $\;\dim H\,$.}

\pf Ecrivons $\;A=\dst\bigoplus_{i\in I} M_{n_i}(\C)\,$. Pour $\;i\in
I\;$ notons $\;r_i\;$ la multiplicit\'e de la repr\'esentation
correspondante. On a alors $\;\dim H=\dst\sum_i n_ir_i\,,\;\dim
A=\dst\sum_i n_i^2\,,\;\dim A'=\dst\sum_i r_i^2\,$. Par
Cauchy-Schwartz, $\;\dim A\dim A'\ge (\dim H)^2\;$ avec \'egalit\'e si
et seulement s'il existe $\;\lambda \in \R\;$ tel que, pour tout $\;i\,$,
$\;r_i=\lambda n_i\,$. Remarquons que $\;H^k\;$ est un $\;A\;$ module
libre si et seulement s'il existe $\;\ell\in\N\;$ tel que, pour tout
$\;i\in I\,$, $\;kr_i=\ell n_i\,$. La premi\`ere assertion en r\'esulte.

Dans ce cas, $\;H\;$ lui m\^eme est libre sur $\;A\,$, si et seulement si
$\;\lambda=r_i/n_i\;$ est entier. Or $\;\dim H=\dst\sum_i
n_ir_i=\lambda\dst\sum_i n_i^2=\lambda \dim A\,$.

\cqfd

Le th\'eor\`eme r\'esulte imm\'ediatement des lemmes \sn 7 et \sn10 (en
utilisant le lemme \sn 6 et la prop. 3.3.d)\ref.

\cqfd

\cor{Soient $\;f\;$ et $\;g\;$ des pr\'e-sous-groupes tels que $\;f\prec
g\,$. Alors $\;D_f\;$ est un $\;D_g\;$ module libre.}

\pf En effet, par le th\'eor\`eme \sn 9, $\;\H\;$ est isomorphe
$\;D_f^k\;$ en tant que $\,D_f$-module, donc en tant que
$\,D_g$-module. Donc $\;D_f^k\;$ est un $\,D_g$-module libre. Le
corollaire r\'esulte alors de la prop. 3.9.c\ref).

\cqfd

\cor{Soient $\;f\,,\;f'\,,\;g\;$ et $\;g'\;$ des pr\'e-sous-groupes, tels
que $\;f\prec f'\;$ et $\;g'\prec g\,$. Alors $\;B_{f,g}\;$ est un
$\;B_{f',g'}$-module libre.}

\pf Notons $\;k=\dim B_{f,g}/\dim B_{f',g}=\dim D_f/ \dim
D_{f'}\,$. Par le corollaire \sn 11, il existe $\;x_1\,...,x_k\;$ tels que
$\;D_f\;$ soit engendr\'e par les $\;x_i\;$ en tant que
$\;D_{f'}\;$ module. Alors $\;B_{f,g}\;$ est engendr\'e par les
$\;x_i\;$ en tant que $\;B_{f',g}$-module. Par un argument de dimension,
il en r\'esulte que $\;B_{f,g}\;$ est un $\;B_{f',g}$-module libre.
Rempla\c cant $\;V\;$ par $\;\Sigma V^*\Sigma\,$, on en d\'eduit que
$\;B_{f',g}\;$ est un $\;B_{f',g'}$-module libre, d'o\`u le r\'esultat.

\cqfd

En particulier, $\;{\cal L}(\H)=B_{\hat e,e}\;$ est un $\;B_{f,g}$-module
libre. Donc,  $\;\H^k\;$ est un $\;B_{f,g}$-module libre, pour tout
$\;k\in\N\;$ tel que la dimension de $\;B_{f,g}\;$ divise $\;kn\;$ \note{Ce
r\'esultat se d\'eduit directement des lemmes \sn 7 et \sn 10. Notons qu'on
peut aussi en d\'eduire une autre d\'emonstration du corollaire \sn 12
(analogue \`a celle de \sn 11).}. Il s'ensuit que $\;\H\;$ est un
$\;B_{f,f}$-module libre de rang 1. Ce r\'esultat est analogue \`a [\ref\Ma].
On peut en fait donner ici une d\'emonstration tr\`es simple de [\ref\Ma].
Par le lemme 3.2.b\ref), les alg\`ebres $\;G_f\;$ et $\;\widehat D_f\;$
commutent.

\prop{\cf [\Ma] th. 3.1\ref) Le $\;G_f\otimes \widehat D_f$-module $\;\H\;$
est isomorphe \`a $\;G_f\otimes \widehat D_f\,$.}

\pf Par la prop. 1.\ref 7.a), la repr\'esentation $\;\pi :G_f\otimes
\widehat D_f\ra \L(\H)\;$ a m\^eme trace normalis\'ee que la
repr\'esentation naturelle $\;\pi ':G_f\otimes \widehat D_f\ra \L(\H\otimes
\H)\,$.  Par le th\'eor\`eme \sn 9, celle-ci co\"incide avec la trace
normalis\'ee de l'action de $\;G_f\otimes \widehat D_f\;$ dans lui m\^eme.
Or  par la prop. 3.3.\ref d), $\;\H\;$ et $\;G_f\otimes \widehat D_f\,$ ont
m\^eme dimension.

\cqfd

\sec{Sous-groupes, co-sous-groupes et sous-quotients.}

\df{ Un pr\'e-sous-groupe $\;f\;$ est appel\'e {\it sous-groupe}\/
(resp. {\it co-sous-groupe}\/) si $\;L_f\;$ (resp. $\;\rho _f\,$) est
dans le centre de $\;S\;$ (resp. $\;\widehat S\,$); il est appel\'e {\it
sous-groupe normal}\/ s'il est \`a la fois un sous-groupe et un
co-sous-groupe. On appelle {\it sous-quotient}\/ de $\;V\;$ un
sous-espace non nul $\;H\;$ de $\;{\cal H}\;$ tel que $\;V(H\otimes
H)=H\otimes H\,$.}

Par la prop. 3.11\ref, un sous-groupe est donc donn\'e par un
projecteur central $\;p\in S\;$ tel que $\;p\otimes p\le \delta (p)\,$.
Remarquons que par la proposition 3.3.c)\ref, on a alors $\;\kappa
(p)=p\,$; on en d\'eduit que notre notion de sous-groupe co\"incide avec
celle de [\Kac].

\prop{Soit $\;H\;$ un sous-quotient.\hfill\break
 a) Il existe un unique pr\'e-sous-groupe $\;f\in H\;$ et un unique
pr\'e-sous-groupe $\;\hat f\in H\;$ tels que $\;H\subset H^f\;$ et
$\;H\subset H_{\hat f}\,$.\hfill\break
 b) On a $\;H=H^f\cap H_{\hat f}\,$.}

\pf a) L'unitaire $\;V\;$ restreint \`a $\;H\otimes H\;$ est
multiplicatif et il admet donc un vecteur fixe \ie un vecteur $\;f\in
H\,,$ non nul tel que $\;V(f\otimes \xi )=f\otimes \xi \;$ pour tout $\;\xi
\in H\,$. Comme $\;V(f\otimes f)=f \otimes f\,$, on peut quitte \`a
multiplier $\;f\;$ par un scalaire, supposer que $\;f\;$ est un
pr\'e-sous-groupe. On a clairement $\;H\subset H^f\,$. Si $\;g\;$ est un
pr\'e-sous-groupe v\'erifiant la m\^eme propri\'et\'e on a $\;f\prec g\;$ et
$\;g\prec f\;$ donc $\;f=g\;$ (3.7.b\ref). L'existence et l'unicit\'e de $\;\hat
f\;$ se d\'emontrent de fa\c con analogue.

\noindent b) Il est clair que $\;H\subset H^f\cap H_{\hat f}\,$. Or d'apr\`es la
prop. 3.9.c\ref), on a $\;\dim(H^f\cap H_{\hat f})\langle f,\hat f\rangle
^2=1\,$; comme l'unitaire $\;V\;$ restreint \`a $\;H\;$ est multiplicatif et
que l'espace des vecteurs fixes est r\'eduit \`a $\;\C f\,$, on a
$\;\dim(H)\langle f,\hat f\rangle ^2=1\;$ (\cf 1.2\ref).

\cqfd

\prop{Soient $\;f,\hat f\;$ des pr\'e-sous-groupes tels que $\;\hat
f\prec f\,$; posons $\;H=H^f\cap H_{\hat f}\,$. Les conditions
suivantes sont
\'equivalentes:\hfill\break
 (i) $\;H\;$ est un sous-quotient.\quad
 (ii) $\;V(H\otimes H^f)= H\otimes H^f\,$.\quad
 (iii) $\;V(H_{\hat f}\otimes H)=H_{\hat f}\otimes H\,$.\quad
 (iv)  $\;V(H_{\hat f}\otimes H^f)= H_{\hat f}\otimes H^f\,$.\quad
 (v) $\;V({\hat f}\otimes f)\in H_{\hat f}\otimes H^f\,$.\quad
 (vi) $\;UH=H\,$.}

\pf Si (ii) est vrai, comme $\;V(H\otimes H)\subset V(H_{\hat
f}\otimes H_{\hat f})\subset {\cal H}\otimes H_{\hat f}\,,$ (lemme
3.1.b), on a
$\;V(H\otimes H)\subset ({\cal H}\otimes H_{\hat f})\cap H\otimes
H^f=H\otimes H\,$. Donc (ii)$\Ra$(i).

Si (iv) est vrai, comme $\;V(H_{\hat f}\otimes H)\subset V(H_{\hat
f}\otimes H_{\hat f})\subset {\cal H}\otimes H_{\hat f}\,,$ (lemme
3.1.b), on a $\;V(H_{\hat f}\otimes H)\subset ({\cal H}\otimes H_{\hat
f})\cap H_{\hat f}\otimes H^f=H_{\hat f}\otimes H\,$. Donc
(iv)$\Ra$(iii).

En rempla\c cant $\;V\;$ par $\;\Sigma V^*\Sigma\,$, on en d\'eduit
que (iv)$\Ra $(ii) et (iii)$\Ra $(i) .

Supposons que (v) soit v\'erifi\'ee. Alors $\;V(\rho _{\hat f}\otimes
L_f) (\hat e\otimes e)\in H_{\hat f}\otimes H^f\,$. Comme $\;\hat e
\otimes e\;$ est s\'eparateur pour $\;\widehat S\otimes S\,$, il en
r\'esulte que
$\;(\rho _{\hat f}\otimes L_f)V(\rho _{\hat f}\otimes L_f)=V(\rho
_{\hat f}\otimes L_f)\,$. Donc (v)$\Ra$(iv).

Montrons que (vi)$\Ra$(v). Par le lemme 3.1.b), on a $\;V(\hat f\otimes
f)\in \H\otimes H_{\hat f}\;$ et $\;\widehat V(f\otimes \hat f)\in
\H\otimes H^f\,$, d'o\`u l'on d\'eduit que $\;V(\hat f\otimes f)\in
UH^f\otimes H_{\hat f}\,$. Comme de plus $\;V\;$ commute avec
$\;\lambda _{\hat f}\otimes R_f\,$, on a $\;V(\hat f\otimes f)\in
(UH_{\hat f}\cap UH^f)\otimes UH^f\cap H_{\hat f}\,$.

Supposons que (vi) est v\'erifi\'ee. Comme $\;R_f\rho_{\hat f}R_f\;$
est l'identit\'e sur $\;H\;$ (rappelons que $\;R_f=UL_fU\,$), on a
$\;R_f\rho_{\hat f}R_f\ge L_f\rho_{\hat f}\,$; or
$\;\tau(R_f\rho_{\hat f}R_f)=\tau(R_f)\tau(\rho_{\hat f})\;$ (prop.
\ref 1.8 appliqu\'ee \`a $\;\widetilde V\,$); il vient $\;L_f\rho _{\hat
f}=R_f\rho _{\hat f}\,$, donc $\;UH^f\cap H_{\hat f}=H\,$; donc
$\;V(\hat f\otimes f)\in H\otimes H\,$.

Supposons enfin que (i) est v\'erifi\'ee. En appliquant 1.3 \`a la
restriction de $\;V\;$ \`a $\;H\otimes H\,$, on voit que pour $\;\xi\in
H\;$ on a $\;(\id\otimes \omega_{\hat f,f})(\Sigma V)\xi =(\dim
H^{-1/2})\xi\,$. D'autre part, $\;(\id\otimes \omega_{\hat f,f})(\Sigma
V)=\langle e,f\rangle^{-1} \langle \hat e,\hat f\rangle^{-1} (\id\otimes
\omega_{\lambda_{\hat f}(\hat e),R_fe})(\Sigma V)= \langle e,f
\rangle^{-1} \langle \hat e,\hat f\rangle^{-1}
R_f(\id\otimes \omega_{\hat e,e})(\Sigma V)\lambda_{\hat
f}= \langle e,f \rangle^{-1} \langle \hat e,\hat f\rangle^{-1}\langle
e,\hat e\rangle R_f\lambda _{\hat f}\,$. Donc $\;R_f\lambda_{\hat
f}\;$ agit comme l'identit\'e sur $\;H\,$. Donc (i)$\Ra$(vi).

\cqfd

Soit $\;H=H^f\cap H_{\hat f}\;$ un sous-quotient. Remarquons que,
$\;(\omega_{f,\hat f}\otimes \id)(\Sigma V)=\langle e,f \rangle^{-1}
\langle \hat e,\hat f\rangle^{-1}(\omega_{R_fe,\lambda _{\hat f}\hat
e}\otimes \id)(\Sigma V)=\langle e,f \rangle^{-1} \langle \hat e,\hat
f\rangle^{-1}(\omega_{e,\hat e}\otimes \id)((R_f\otimes 1)\Sigma
V(\lambda _{\hat f}\otimes 1))=\langle e,f \rangle^{-1} \langle \hat
e,\hat f\rangle^{-1}(\omega_{e,\hat e}\otimes \id)((1\otimes \lambda
_{\hat f})\Sigma V(1\otimes R_f))=\langle e,f \rangle^{-1} \langle
\hat e,\hat f\rangle^{-1}\langle e,\hat e \rangle \lambda _{\hat
f}UR_f\;$ par 1.3. On d\'eduit alors de 1.3\ref\ que l'op\'erateur $\;\og
U\fg\;$ associ\'e \`a la restriction de
$\;V\;$ \`a $\;H\,$ est la restriction de $\;U\;$ \`a $\;H\,$.

\cor{Soient $\;f\;$ et $\;g\;$ des pr\'e-sous-groupes tels que
$\;f\prec g\;$ et $\;\langle f,g\rangle=2^{-1/2}\,$. Alors $\;H_f\cap
H^g\;$ est un sous-quotient.}

\pf On a $\;Uf=f\;$ et $\;Ug=g\,$, donc $\;U(H_f\cap H^g)=H_f\cap
H^g\,$, vu que cet espace est de dimension 2. Le r\'esultat d\'ecoule de
la prop. \sn3.

\cqfd

\cor{Soient $\;H\;$ un sous-quotient, $\;f,\hat f\in H\;$ les
pr\'e-sous-groupes tels que $\;H=H^f\cap H_{\hat f}\,$.\hfill\break
 a) Si $\;\hat e\in H\;$ on a $\;\hat f=\hat e\,,\; H^f=H\;$ et $\;V({\cal
H}\otimes H)={\cal H}\otimes H\,$. Le projecteur $\;q\;$ sur $\;H\;$
est dans le centre de $\;S\,$. \hfill\break
 b) Si $\;e\in H\;$ on a $\;f=e\,,\; H_{\hat f}=H\;$ et $\;V(H\otimes
{\cal H})=H\otimes {\cal H}\,$. Le projecteur sur $\;H\;$ est dans le
centre de $\;\widehat S\,$.}

\pf a) L'\'egalit\'e $\;\hat f=\hat e\;$ r\'esulte de l'unicit\'e de
$\;\hat f\,$; comme $\;H_{\hat e}=\H\,$, on a $\;H^f=H\;$ et
l'\'egalit\'e $\;V({\cal H}\otimes H)={\cal H}\otimes H\;$
(\'equivalente \`a $\;q\in S'\,$) r\'esulte de \sn3. Comme $\;q=L_f\,,\;
q\in S\,$; donc $\;q\in S\cap S'\,$.

\noindent b) se d\'eduit de a) en rempla\c cant $\;V\;$ par $\;\Sigma
V^*\Sigma \,$.

\cqfd

\cor{Pour un pr\'e-sous-groupe $\;f\;$ les propri\'et\'es suivantes sont
\'equivalentes.\hfill\break
 (i) $\;f\;$ est un co-sous-groupe.\hfill\break
 (ii) $\;V(H_f\otimes \H)=H_f\otimes \H\,$.\hfill\break
 (iii) $\;V(H_f\otimes H_f)=H_f\otimes H_f\,$.\hfill\break
 (iv) $\;\sigma(\delta (L_f))=\delta (L_f)\,$.\hfill\break
 (v) $\;UH_f=H_f\,$.\hfill\break
 (vi) $\;D_f=G_f\,$.}

\pf Comme la transform\'ee de Fourier de $\;\rho _f\;$ est
proportionnelle \`a $\;L_f\,$, l'\'equivalence entre (i) et (iv)
r\'esulte de la prop. 1.5\ref. (i)$\iff$(ii) est claire. Enfin (ii)$\iff$(iii)
r\'esulte imm\'ediatement de la prop. \sn3. Si (v) est satisfaite,
$\;\rho _f=\lambda _f\,$, donc (v)$\Ra$(i); par la prop. \sn3, tout
sous-quotient est invariant par $\;U\;$ donc (ii)$\Ra$(v). Enfin
(v)$\iff$(vi) r\'esulte des d\'efinitions de $\;G_f\;$ et $\;D_f\,$.

\cqfd

En particulier (prop. 3.11), les co-sous-groupes de $\;V\;$ sont donn\'es
par les projecteurs $\;p\in S\;$ tels que $\;\sigma(\delta (p))=\delta
(p)\ge p\otimes p\,$. On retrouve la notion de `factor group' de [\Kac].
Donc un sous-groupe normal est donn\'e par un projecteur central
$\;p\in S\;$ tel que $\;\sigma(\delta (p))=\delta (p)\ge p\otimes p\,$.
Cette notion co\"incide donc avec celle de [\Kac].

On a aussi \'evidemment un analogue du corollaire \sn6 pour les
sous-groupes.

\prop{La borne sup\'erieure et la borne inf\'erieure de deux
sous-groupes (resp. co-sous-groupes) est un sous-groupe (resp.
co-sous-groupe).}

\pf Soient $\;f\;$ et $\;g\;$ deux sous-groupes (resp. co-sous-groupes).
Notons $\;h\;$ leur borne inf\'erieure. Notons $\;A\;$ le centre de
$\;S\;$ (resp.  $\;A=\{\,x\in S\,,\;\delta (x)=\sigma(\delta (x))\,\}\,$).
Par hypoth\`ese, $\;L_f\in A\;$ et $\;L_g\in A\,$, donc le projecteur
orthogonal sur $\;H^h=H^f\cap H^g\;$ est un \'el\'ement de
$\;A\,$, donc $\;h\;$ est un sous-groupe (resp. un co-sous-groupe).
Quand on change $\;V\;$ en $\;\Sigma V^*\Sigma\,$, les
pr\'e-sous-groupes ne changent pas, l'ordre est chang\'e en son
oppos\'e, les sous-groupes deviennent des co-sous-groupes {\it et
vice-versa,}\/ d'o\`u la deuxi\`eme assertion.

\cqfd

\ssec{Normalisateur d'un pr\'e-sous-groupe.}

Soit $\;f\;$ un pr\'e-sous-groupe.

\lem{Pour $\;\xi \in H^f\;$ et $\;\zeta \in {\cal H}\,,$ on a $\; V(\xi
\otimes \zeta )\in {\cal H}\otimes H^f\iff V(\xi \otimes \zeta )\in
H^f\otimes {\cal H}\,$.}

\pf Par la prop. 3.11.a\ref), on a $\;(L_f\otimes 1)V(L_f\otimes
1)=(1\otimes L_f)V(L_f\otimes 1)=(L_f\otimes L_f)V\,$. Donc
$\;(L_f\otimes 1)V(\xi \otimes \zeta )=V(\xi \otimes \zeta )\iff
(1\otimes L_f)V(\xi \otimes \zeta )=V(\xi \otimes \zeta )\,$.

\cqfd

\prop{Soit $\;f\;$ un pr\'e-sous-groupe. Posons $\;H = \{\, \xi \in
H^f\,,\; V(\xi \otimes H^f)\subset {\cal H}\otimes H^f \,\}\,$.
\hfill\break
 a) On a $\;V(H\otimes H)=H\otimes H\,$, autrement dit $\;H\;$ est un
sous-quotient.\hfill\break
 b) $\;H\;$ est le plus grand sous-quotient de vecteur fixe $\;f\;$ \ie
tel que $\;f\in H\subset H^f\,$.\hfill\break
 c) L'ensemble $\;K = \{\, \xi \in H_f\,,\; V^*(H_f\otimes \xi)\subset
H_f\otimes  \H\,\}\,$  est le plus grand sous-quotient de vecteur
cofixe $\;f\;$ \ie tel que $\;f\in K\subset H_f\,$}

\pf a) On a $\;V(H\otimes H^f)\subset {\cal H}\otimes H^f\;$ par
d\'efinition de
$\;H\,$. On d\'eduit alors du lemme \sn8 que $\;V(H\otimes H^f)\subset
H^f\otimes {\cal H}\,$.

Pour $\;\xi ,\zeta \in H\;$ et $\;\alpha \in H^f , V_{23}(V(\xi \otimes
\zeta )\otimes \alpha ) = V_{12}V_{13}V_{23}(\xi
\otimes \zeta \otimes \alpha ) \in {\cal H}\otimes {\cal H}\otimes
H^f\,$, donc $\;V(\xi \otimes \zeta )\in {\cal H}\otimes H\,$.

Pour $\;\xi \in H\;$ et $\;\alpha ,\zeta \in H^f ,\; V_{12}V_{13}(\xi
\otimes \alpha \otimes \zeta )=V_{23}V_{12}(\xi \otimes V^*(\alpha
\otimes \zeta )) \in H^f\otimes {\cal H}\otimes {\cal H}\;$ car
$\;V^*(H^f\otimes H^f) \subset H^f\otimes {\cal H}\;$ (lemme 3.1.b\ref); il
r\'esulte du lemme \sn 8 que $\;V(\xi \otimes \zeta )\in H\otimes {\cal
H}\,$.

\medskip\noindent b) Si $\;H_1\;$ est un autre tel sous-quotient, on a
$\;V(H_1\otimes H^f)=H_1\otimes H^f\,$, (prop. \sn 3) donc
$\;H_1\subset H\,$.

\medskip\noindent c) se d\'eduit de a) et b) en rempla\c cant $\;V\;$
par $\;\Sigma V^*\Sigma\,$.

\cqfd

Le pr\'e-sous-groupe $\;g\;$ vecteur fixe (resp. co-fixe) du plus grand
sous-quotient de vecteur co-fixe (resp. fixe) $\;f\;$ s'appelle le {\it
nomalisateur }\/ (resp {\it co-normalisateur)}\/ de $\;f\,$.

\cor{Soit $\;\hat f\;$ un pr\'e-sous-groupe. L'ensemble des
pr\'e-sous-groupes $\;f\;$ satisfaisant $\;\hat f\prec f\;$ et tels que
$\;H_{\hat f}\cap H^f\;$ soit un sous-quotient est stable par borne
sup\'erieure et par borne inf\'erieure.}

\pf Soit $\;H\;$ le plus grand sous-quotient de vecteur cofixe $\;\hat
f\,$. L'ensemble, ordonn\'e par l'inclusion, des sous-quotients de
vecteur cofixe $\;\hat f\;$ s'identifie avec l'ensemble des
sous-groupes de la restriction de $\;V\;$ \`a $\;H\otimes H\,$. Le
r\'esultat d\'ecoule alors ais\'ement de la prop. \sn 7.

\cqfd

\sec {Une g\'en\'eralisation du \og biproduit crois\'e\fg.}

Dans ce paragraphe, on fixe un unitaire multiplicatif $\;V\;$ de
multiplicit\'e  $\;1\,$ dans un espace hilbertien $\;\H\;$ de dimension finie
$\;n\;$ et deux pr\'e-sous-groupes $\;f\;$ et $\;g\;$ de $\;V\;$ tels que
$\;\langle f,g\rangle =n^{-1/2}\,$. Le principal r\'esultat est que les
alg\`ebres $\;A_{f,f}\;$ et $\;B_{g,g}\;$ introduites au paragraphe 4\ref\
sont des alg\`ebres de Hopf en dualit\'e (th\'eor\`eme \sn2). Cette
construction g\'en\'eralise la construction du \og biproduit crois\'e\fg\ de
[\Kac], [\Takeu], [\Maj]  (voir aussi [\BSb] et [\HS]).

Commen\c cons par quelques propri\'et\'es \'el\'ementaires v\'erifi\'ees
par le couple $\;(f,g)\;$ de pr\'e-sous-groupes tels que $\;\langle f,g\rangle
=n^{-1/2}\,$; certaines d'entre elles nous serviront dans la d\'emonstration
du th\'eor\`eme \sn 2.

\prop{On a les identit\'es suivantes\hfill\break
 a) $\;L_fL_g=L_gL_f=\theta _{\hat e,\hat e}\;$ et $\;\rho _f\rho _g=\rho
_g\rho _f=\theta _{e,e}\,$.\hfill\break
 b) $\;\langle f,e\rangle =\langle g,\hat e\rangle \;$ et $\;\langle
f,\hat e\rangle =\langle g,e\rangle \,$.\hfill\break
 c) $\;\rho _f(g) =\langle g,e\rangle e\;$ et $\;L_g(f)=\langle f,\hat e\rangle
\hat e\,$.\hfill\break
 d) $\;\rho _fL_g\rho _f=\langle g,e\rangle ^2\rho _f\;$ et $\;L_g\rho
_fL_g=\langle f,\hat e\rangle ^2L_g\,$.\hfill\break
 e) $\;L_gR_f=\theta _{\hat e,\hat e}\;$ et $\;\lambda _f\rho _g=\theta
_{e,e}\,$.}

\pf a) Posons $\;x=L_fL_g-\theta _{\hat e,\hat e}=L_f(L_g-\theta _{\hat
e,\hat e})\;$ et $\;y=\rho _f\rho _g-\theta _{e,e}\,$. On a
$\;x^*x=L_gL_fL_g-\theta _{\hat e,\hat e}\;$ et $\;y^*y=\rho _g\rho _f\rho
_g-\theta _{e,e}\,$. De plus $\;\tau (x^*x)=\tau (L_fL_g)-1/n=\langle
L_fe,L_ge\rangle -1/n=\langle e,f\rangle \langle g,e\rangle \langle
f,g\rangle -1/n\,$. De m\^eme $\;\tau (y^*y)=\langle \hat e,f\rangle \langle
g,\hat e\rangle \langle f,g\rangle -1/n\,$. Or $\;\langle e,f\rangle \langle
g,e\rangle \langle f,g\rangle\langle \hat e,f\rangle \langle
g,\hat e\rangle \langle f,g\rangle=1/n^2\,$. On en d\'eduit que $\;\tau
(x^*x)=\tau (y^*y)=0\,$, d'o\`u a) et $\;\langle \hat e,f\rangle \langle
g,\hat e\rangle \langle f,g\rangle=1/n\,$, d'o\`u b).

\noindent c) On a $\;\rho _f(g)=\rho _f\rho _g(g)\;$ et
$\;L_g(f)=L_gL_f(f)\,$. c) d\'ecoule donc de a).

\noindent d) On a
$$\eqalign
{\rho _fL_g\rho _f&=(\omega _{g,g}\otimes
\id\otimes \omega _{f,f}\otimes \omega _{f,f})(V_{23}V_{12}V_{24})\cr
& =(\omega _{g,g}\otimes
\id\otimes \omega _{f,f}\otimes \omega
_{f,f})(V_{12}V_{13}V_{23}V_{24})\cr
& =(\omega _{g,g}\otimes \id\otimes \omega _{f,f}\otimes \omega
_{f,f})(V_{12}V_{13}V_{34}V_{23})\cr
& =(\omega _{g,g}\otimes \id\otimes \omega
_{f,f})(V_{12}V_{13}(1\otimes 1\otimes \rho _f)V_{23}) \,.}$$
Or $$\eqalign
{(\id\otimes\id\otimes \omega _{f,f})(V_{13}(1\otimes 1\otimes \rho
_f)V_{23})&= (\id\otimes\id\otimes \omega _{f,f})(V_{13}(1\otimes
1\otimes \rho _f)V_{23}(1\otimes 1\otimes R _f))\cr&=
(\id\otimes\id\otimes \omega _{f,f})(V_{13}(1\otimes 1\otimes \theta
_{f,f})V_{23})\cr&=\rho _f\otimes\rho _f\,.}$$ Donc $\;\rho _fL_g\rho
_f=(\omega _{g,g}\otimes \id)(V(\rho _f\otimes \rho _f))=(\omega _{g,\rho
_fg}\otimes \id)(V)\rho _f=\langle g,e\rangle (\omega _{g,e}\otimes \id
)(V)\rho _f\,$.

La deuxi\`eme assertion s'en d\'eduit en rempla\c cant $\;V\;$ par
$\;\Sigma V^*\Sigma\,$.

\noindent e) $\;L_gR_f\;$ est un projecteur plus grand que $\;\theta _{\hat
e,\hat e}\,$; par ailleurs $\;L_gR_f\;$ est proportionnel \`a $\;L_g\rho
_fL_gR_f=L_g\rho _fR_fL_g=L_g\theta _{f,f}L_g\;$ qui est de rang $\;1\,$,
d'o\`u l'\'egalit\'e. La deuxi\`eme assertion s'en d\'eduit en rempla\c cant
$\;V\;$ par $\;\widehat V\,$.

\cqfd

Venons au r\'esultat principal de ce paragraphe:

\th{Il existe un (et un seul) unitaire multiplicatif $\;W\in \L(\H\otimes
\H)\;$ de multiplicit\'e $\;1\,$, pour lequel $\;g\;$ est fixe et $\;f\;$ cofixe
et tel que les $\,C^*$-alg\`ebres $\;S\;$ et $\;\widehat S\;$ associ\'ees sont
respectivement $\;A_{f,f}\;$ et $\;B_{g,g}\,$.}

Nous avons not\'e qu'un unitaire multiplicatif dans $\;\L(\H\otimes \H)\;$ de
multiplicit\'e $\;1\;$ est d\'etermin\'e par l'espace de ses vecteurs fixes,
de ses vecteurs cofixes et les $\,C^*$-alg\`ebres $\;S\;$ et $\;\widehat
S\;$ associ\'ees (\cf prop. 1.10\ref); donc si un tel op\'erateur $\;W\;$
existe, il est unique.

Notons $\;\tau \;$ la trace normalis\'ee de $\;\L(\H)\,$. Le th\'eor\`eme
\sn2 r\'esulte de la caract\'erisation des alg\`ebres de Hopf en dualit\'e
donn\'ee au paragraphe 2 (th\'eor\`eme 2.8\ref). On a $\;\theta
_{f,f}=L_f\lambda _f\in A_{f,f}\;$ et $\;\theta _{g,g}=L_g\rho _g\in
B_{g,g}\,$. De plus, $\;\dim A_{f,f}=\dim B_{g,g}=n\;$ (lemme 4.6,
prop. 3.3.d\ref). Nous devons donc montrer:

\noindent --- Pour tout $\;a\in A_{f,f}\;$ et tout $\;b\in B_{g,g}\;$
on a $\;\tau (ab)=\tau (a)\tau (b)\;$ (prop. \sn \ref4 ci-dessous).

Il en r\'esulte qu'il existe un unitaire $\;W\in\L(\H\otimes \H)\;$ tel que
$\;\langle ag\otimes bf,W(cg\otimes df)\rangle =\tau (a^*b^*cd)\;$ pour tout
$\;a,c\in A_{f,f}\;$ et tout $\;b,d\in B_{g,g}\,$.

\noindent --- On a $\;W\in B_{g,g}\otimes \L(\H)\,$. En vertu de la
prop. 4.8\ref, il revient au m\^eme de d\'emontrer que, pour tout $\;T\in
R(D_g)\,$, on a $\;[T\otimes 1,W]=0\,$, (prop. \sn\ref7) et que  pour tout
$\;T\in \lambda (\widehat G_g)\,$, on a $\;[T\otimes 1,W]=0\;$ (prop. \sn
11).

\lem{Pour tout $\;a\in  \lambda(\widehat {D_f})\,$, tout $\;x\in S\;$ et
tout $\;b\in \widehat {G_g}\;$ on a $\;\tau (axb)=\tau (a)\tau(x)\tau (b)\,$.}

\pf Consid\'erons les sous-alg\`ebres $\;D_1= \lambda(\widehat
{D_f})\otimes 1\,,\,D_2 = L(S)\otimes 1\,,\,D_3 =
V^*(1\otimes\rho(\widehat {G_g}))V\,$.\hfill\break
Remarquons que l'espace vectoriel engendr\'e par $\;D_1D_2\;$ est
$\;A_{f,\hat e}\otimes 1\;$ et que l'espace vectoriel engendr\'e par
$\;D_2D_3\;$ est $\;\widetilde V(B_{\hat e,g}\otimes 1)\widetilde V^*\,$;
de plus, $\;D_1\;$ et $\;D_3\;$ commutent; donc l'espace engendr\'e par
$\;d_1d_2d_3\,,\;d_i\in D_i\;$ est une sous-$C^*$-alg\`ebre $\;D\;$ de
$\;\L(\H)\otimes \widehat S\,$. La dimension de $\;D\;$ est clairement
$\;\le n^2\,$.

L'application $\;\sigma =\id\otimes \varepsilon :D\ra \L(\H)\;$ est une
repr\'esentation. Son image contient $\;\lambda _f\rho_g=\theta_{e,e}\;$ et
$\;S\;$ donc tout $\;\L(\H)\,$, vu que $\;e\;$ est totalisateur pour $\;S\,$.
Pour des raisons de dimension, on en d\'eduit que $\;\sigma \;$ est un
isomorphisme. Comme $\;D\;$ est un facteur, les traces $\;\tau \otimes
\tau \;$ et $\;\tau \circ \sigma \;$ co\"incident.

On a $\,(\tau\otimes \tau)(ax\otimes 1)
V^*(1\otimes b)V)=\tau \circ \sigma ((ax\otimes 1)
V^*(1\otimes b)V)=\tau (axb)\,$. Or $\;\tau \;$ est la mesure
de Haar de  $\;\widehat S\,$, donc $\;(\id\otimes \tau) (V^*(1\otimes
b)V)=\tau (b)1\,$; donc $\;\tau (axb)=\tau (b)\tau
(ax)=\tau (b)\tau (a)\tau (x)\;$ (prop. 1.7.a\ref).

\cqfd

\prop{Pour tout $\;a\in A_{f,f}\;$ et tout $\;b\in B_{g,g}\;$ on a $\;\tau
(ab)=\tau (a)\tau (b)\,$.}

\pf Pour $\;a_1\in
\lambda (\widehat D_f)\,,\; a_2\in L(G_f)\,,\;b_2\in L(D_g)\,,$ et $\;b_3\in
\rho (\widehat G_g)\,$, on a $\;\tau (a_1a_2b_2b_3)=\tau (a_1)\tau
(a_2b_2)\tau (b_3)\;$ (lemme \sn3\ref). Or  $\;\tau (a_2b_2)=\langle
e,a_2b_2e\rangle =\langle \rho _fe,a_2b_2\lambda _ge\rangle\,$; or $\;\rho
_fa_2=a_2\rho _f\;$ et $\;b_2\lambda _g=\lambda _gb_2\;$ (prop. 4.2\ref);
enfin $\;\rho _f\lambda _g=\theta _{e,e}\,$; on trouve $\;\langle \rho
_fe,a_2b_2\lambda _ge\rangle=\langle e,a_2e\rangle \langle e,b_2e\rangle
\,$. On a donc $\;\tau (a_1a_2b_2b_3)=\tau (a_1)\tau (a_2)\tau (b_2)\tau
(b_3)=\tau (a_1a_2)\tau (b_2b_3)\,$, d'o\`u la proposition.

\cqfd

\lem{Pour tout $\;b\in B_{\hat e,g}\,$, on a $\;E_{A_{f,\hat e}}(b)\in S\,$.}

\pf Pour $\;a\in A_{f,\hat e}\,,\;x\in S\,,\;b\in
\widehat G_g\,$, on a $\;ax\in A_{f,\hat e}\,$, donc $\;\tau (axb)=\tau
(ax)\tau (b)\;$ (lemme \sn 3), donc $\;E_{A_{f,\hat e}}(xb)=\tau (b)x\,$.

\cqfd

 Consid\'erons l'op\'erateur $\;W\in \L(\H\otimes \H)\;$ donn\'e par
$\;\langle ag\otimes bf,W(cg\otimes df)\rangle =\tau (a^*b^*cd)\;$ pour
tout $\;a,c\in A_{f,f}\;$ et tout $\;b,d\in B_{g,g}\,$.

Par le th\'eor\`eme 2.8\ref, il nous reste juste \`a montrer que $\;W\in
B_{g,g}\otimes \L(\H)\,$.

Nous aurons besoin du lemme suivant:

\lem{a) On a $\;(\rho_g\otimes 1)\widehat V(1\otimes \lambda
_g)=\widehat V(\rho_g\otimes \lambda _g)\,$.\hfill\break
 b) Pour tout $\;\alpha ,\beta \in \lambda (\widehat D_f)\;$ on a
$\;[\rho _gR_f\otimes 1,(\beta ^*\otimes 1)\widehat V(\alpha \otimes
\lambda _g)]=0\,$.}

\pf a) On a $\;\Sigma(1\otimes U)\widehat V^*(\rho_g\otimes 1)\widehat
V(1\otimes \lambda _g)(1\otimes U)\Sigma=V^*(1\otimes \rho _g)V(\rho
_g\otimes 1)=\widehat \delta (\rho _g)(\rho _g\otimes 1)=\rho _g\otimes
\rho _g\;$ (prop. 3.11.b\ref), d'o\`u a).

\noindent b) Par a) $\;\rho _g\otimes 1\;$ commute avec $\;\widehat
V(1 \otimes \lambda _g)\,$. Les autres commutations sont \'evidentes.

\cqfd

\prop{Pour tout $\;T\in R(D_g)\,$, on a $\;[T\otimes 1,W]=0\,$.}

\pf Soient $\;x,a\in G_f\,$, $\;y\in B_{g,g}\,$, $\;b\in D_g\,,\;c\in
\widehat G_g\,$, et $\;\alpha ,\beta \in \lambda(\widehat D_f)\,$. On a
$\;[T,\rho_g]=0\,$, donc $\;T\alpha g\in H_g=\lambda(\widehat D_f)g\,$;
soit $\;\alpha _1\in \lambda(\widehat D_f)\;$ tel que $\;T\alpha
g=\alpha _1g\,$. De m\^eme, il existe un unique $\;\beta _1\in \lambda(
\widehat D_f)\;$ tel que $\;T^*\beta g=\beta _1g\,$. On a
$$\eqalign{\langle x\beta g\otimes yf ,(T\otimes 1)W(a\alpha g\otimes
cbf)\rangle&=
\langle x\beta _1g\otimes yf ,W(a\alpha g\otimes cbf)\rangle=\tau (\beta
_1^*x^*y^*a\alpha cb)\cr&=\tau (x^*y^*ac\alpha b\beta_1^*)\,.}$$
De m\^eme
$$\langle x\beta g\otimes yf,W(T\otimes 1)(a\alpha g\otimes cbf)\rangle=
\tau (x^*y^*ac\alpha _1b\beta^*)\,.$$
Remarquons que $\;\alpha b\beta_1^*,\alpha _1b\beta^*\in A_{f,\hat e}\;$
et $\;x^*y^*ac\in B_{\hat e,g}\,$. On a alors  $\;\tau
(x^*y^*ac\alpha b\beta_1^*)=\tau (z\alpha b\beta _1^*)\;$ et $\;\tau
(x^*y^*ac\alpha _1b\beta^*)=\tau(z\alpha _1b\beta^*)\,$, o\`u
$\;z=E_{A_{f,\hat e}}(x^*y^*ac)\,$. Par le lemme \sn5\ref, $\;z\in S\,$.

Comme $\;z,b\in S\;$ et $\;\alpha _1,\beta ^*\in \lambda (\widehat S)\,$, on
a $\;\tau (\beta^*z\alpha _1b)=\langle \beta \hat e\otimes z^*e,\widehat
V(\alpha _1\hat e\otimes be)\rangle \,$; de m\^eme, $\;\tau
(\beta_1^*z\alpha b)=\langle \beta _1\hat e\otimes z^*e,\widehat
V(\alpha \hat e\otimes be)\rangle \,$. Or $\;be\in UH_g\;$ et $\;\hat
e=\langle g,\hat e\rangle ^{-1}R_fg=\langle g,\hat e\rangle ^{-1}R_f\rho
_gg\,$.

Donc $\;\tau (\beta^*z\alpha _1b)=\langle g,\hat e\rangle ^{-1} \langle
g\otimes z^*e,(\rho _gR_f\otimes 1)(\beta ^*\otimes 1)\widehat V(\alpha
_1\otimes \lambda _g)(\hat e\otimes be)\rangle = \langle g,\hat e\rangle
^{-1} \langle g\otimes z^*e,(\beta ^*\otimes
1)\widehat V(\alpha _1\otimes \lambda _g)(\rho _gR_f \hat e\otimes
be)\rangle \;$ par le lemme \sn6\ref.
Donc $\;\tau (\beta^*z\alpha _1b)=
\langle g,\hat e\rangle ^{-2} \langle \beta g\otimes z^*e,\widehat V(\alpha
_1g\otimes be)\rangle\,$. De m\^eme, $\;\tau
(\beta_1^*z\alpha b)= \langle g,\hat e\rangle ^{-2} \langle \beta
_1g\otimes z^*e,\widehat V(\alpha g\otimes be)\rangle= \langle g,\hat
e\rangle ^{-2} \langle T^*\beta g\otimes z^*e,\widehat V(\alpha g\otimes
be)\rangle= \tau (\beta^*z\alpha _1b)\,$, d'o\`u le r\'esultat.

\cqfd

\lem{Soient $\;a_1\in \lambda(\widehat {D_f})\,,\;a_2\in
L(G_f)\,,\;b_2\in L(D_g)\;$ et $\;b_1\in \rho(\widehat {G_g})\,$.\hfill\break
a) On a $\;\langle a_1a_2g,b_1b_2f\rangle = n^{1/2}\langle a_2e,b_1\hat
e\rangle \langle a_1\hat e,b_2e\rangle\,$.\hfill\break
b) On a $\;\langle Ub_2b_1f,a_2a_1g\rangle = n^{1/2}\langle Ub_1\hat
e,a_2e\rangle \langle Ub_2e,a_1\hat e\rangle\,$.}

\pf a) On a $\;[R_f,a_1]=[R_g,b_1]=0\;$ et $\;R_fR_g=\theta _{\hat e,\hat
e}\,$, donc $$\eqalign{\langle a_1a_2g, b_1b_2f\rangle &= n^{1/2}\langle
b_1^*a_2R_ge,a_1^*b_2R_fe\rangle \cr&
=n^{1/2}\langle R_gb_1^*a_2e,R_fa_1^*b_2e\rangle \cr&
=n^{1/2}\langle b_1^*a_2e, \theta _{\hat e,\hat e}a_1^*b_2e\rangle\cr&
= n^{1/2}\langle a_2e,b_1\hat e\rangle \langle a_1\hat e,b_2e\rangle\,.}$$
 b) On a $\;[\rho _f,a_2]=[\lambda _g,b_2]=0\;$ et $\;\rho
_f\rho _g=\theta _{e,e}\,$, donc $$\eqalign{\langle
Ub_2b_1f,a_2a_1g\rangle &= n^{1/2}\langle
a_2^*Ub_1\lambda _f\hat e,Ub_2^*Ua_1\rho _g\hat e\rangle \cr&
=n^{1/2}\langle \rho _fa_2^*Ub_1\hat e,\rho _gUb_2^*Ua_1\hat e\rangle
\cr&=n^{1/2}\langle a_2^*Ub_1\hat e,
\theta _{e,e}Ub_2^*Ua_1\hat e\rangle \cr&= n^{1/2}\langle Ub_1\hat
e,a_2e\rangle \langle Ub_2e,a_1\hat e\rangle\,.}$$

\cqfd

Notons $\;J_A:xg\mt x^*g\;$ ($\,x\in A_{f,f}\,$) et $\;J_B:yf\mt y^*f\;$
($\,y\in B_{g,g}\,)$ les involutions canoniques.

\prop{On a $\;J_AJ_B=U\,$.}

\pf Soient $\;a_1\in \lambda(\widehat {D_f})\,,\;a_2\in
L(G_f)\,,\;b_2\in L(D_g)\;$ et $\;b_1\in \rho(\widehat {G_g})\,$. On a:
$\;\langle J_AJ_Bb_2b_1f,a_2a_1g\rangle = \langle
J_Aa_2a_1g,J_Bb_2b_1f\rangle= \langle a_1^*a_2^*g,b_1^*b_2^*f\rangle=
n^{1/2}\langle a_2^*e,b_1^*\hat e\rangle \langle a_1^*\hat
e,b_2^*e\rangle\;$ par le a) du lemme pr\'ec\'edent. Or,  $\;\langle
a_2^*e,b_1^*\hat e\rangle=\langle Ja_2e,\widehat Jb_1\hat
e\rangle=\langle J\widehat Jb_1\hat e,a_2e\rangle=\langle Ub_1\hat
e,a_2e\rangle\,$; de m\^eme, $\;\langle a_1^*\hat
e,b_2^*e\rangle=\langle Ub_2e,a_1\hat e\rangle\,.$ Par le lemme
\sn8.\ref b), on a donc $\;\langle J_AJ_Bb_2b_1f,a_2a_1g\rangle =\langle
Ub_2b_1f,a_2a_1g\rangle\,$.

Comme les $\;b_2b_1f\;$ et les $\;a_2a_1g\;$ engendrent $\;\H\,$, le
r\'esultat en d\'ecoule imm\'ediatement.

\cqfd

\lem{Pour $\;a,x\in UA_{f,f}U\;$ et pour $\;b,y\in B_{g,g}\;$ on a: $\;\langle
ag\otimes bf,(1\otimes U)W^*(1\otimes U)(xg\otimes yf)\rangle
=\tau(a^*b^*xy)\,$.}

\pf Posons $\;a_1=J_AaJ_A\,,\;x_1=J_AxJ_A\,,\;b_1=J_AbJ_A \;$
et $\;y_1=J_AyJ_A\,$. Alors $\;a_1,x_1\in A_{f,f}\;$ et
$\;b_1,y_1\in B_{g,g}\,$. On a
$\;\langle ag\otimes bf,(1\otimes U)W^*(1\otimes U)(xg\otimes
yf)\rangle=\langle (J_Aa_1g\otimes
J_Ab_1g),(1\otimes U)W^*(1\otimes U)(J_Ax_1f\otimes
J_Ay_1g)\rangle=\langle (J_Aa_1g\otimes
UJ_Ab_1f),W^*(J_Ax_1g\otimes
UJ_Ay_1f)\rangle=\langle W(a_1^*g\otimes b_1^*f),(x_1^*g\otimes
y_1^*f)\rangle=\tau(b_1a_1y_1^*x_1^*)=\tau (J_Abay^*x^*J_A)=\tau
(xya^*b^*)\,$.

\cqfd

\prop{Pour tout $\;T\in \lambda (\widehat G_g)\,$, on a $\;[T\otimes
1,W]=0\,$.}

\pf Si on remplace $\;V\;$ par $\;\Sigma V^*\Sigma \,,$ alors
les alg\`ebres $\;S\;$ et $\;\widehat S\;$ sont \'echang\'ees, toutes deux
avec le coproduit oppos\'e; les alg\`ebres $\;D_g\;$ et $\;\widehat G_g\;$
sont \'echang\'ees, l'alg\`ebre $\;G_f\;$ est transform\'ee en $\;\widehat
D_f\;$ et $\;U\widehat D_fU\;$ en $\;UG_fU\,$.
Donc $\;B_{g,g}\;$ est inchang\'ee et $\;A_{f,f}\;$ est transform\'ee en
$\;UA_{f,f}U\,$. Par le lemme \sn10, $\;W\;$ est alors transform\'ee en
$\;(1\otimes U)W^*(1\otimes U)\,$.

La prop. \sn7\ref\ montre que $$[(1\otimes U)W^*(1\otimes U)
,(U\widehat G_gU\otimes 1)]=0\,.$$

\cqfd

\rem On peut donner une formule explicite pour l'unitaire multiplicatif
$\;W\,$: notons $\;T:\H\otimes \H\ra \H\;$ l'application donn\'ee par
$\;T(\xi\otimes \eta)=aR_g\rho_f\eta\,$, o\`u $\;a\in U\widehat D_fU\;$
est tel que $\;a\hat e=\langle f,\hat e\rangle ^{-1}R_f\xi\,$. On v\'erifie
que $\;T^*T=R_f\otimes \rho _f\,$, donc $\;T^*\;$ est isom\'etrique. On a
$$W=n(T\otimes T)(1\otimes R_g\otimes \lambda_g\otimes 1)
\widehat V{}_{13}\widehat V{}_{23}V_{23}V_{24} (T^*\otimes T^*)\,.$$

\prop{a) $\;e\;$ est un pr\'e-sous-groupe de
$\;W\,$; on a $\;(\omega_{e,e}\otimes \id)(W)=\lambda_f\;$ et
$\;(\id \otimes \omega_{e,e})(W)=\rho_g\,$. La sous-alg\`ebre co\"id\'eale
\`a gauche (resp. \`a droite) de l'alg\`ebre de Hopf $\;A_{f,f}\;$ (resp.
$\;B_{g,g}\,$) associ\'e est $\;U\widehat D_fU\;$ (resp. $\;\widehat G_g\,$).
\hfill\break
 b) $\;\hat e\;$ est un pr\'e-sous-groupe de
$\;W\,$; on a $\;(\omega_{\hat e,\hat e}\otimes \id)(W)=L_f\;$ et
$\;(\id \otimes \omega_{\hat e,\hat e})(W)=L_g\,$. La sous-alg\`ebre
co\"id\'eale \`a droite (resp. \`a gauche) de l'alg\`ebre de Hopf $\;A_{f,f}\;$
(resp. $\;B_{g,g}\,$) associ\'e est $\;G_f\;$ (resp. $\;D_g\,$).\hfill\break
c) Le pr\'e-sous-groupe $\;e\;$ de $\;W\;$ est un sous-groupe (resp. un
co-cous-groupe) si et seulement si $\;f\;$ (resp. $\;g\,$) est un
co-sous-groupe de $\;V\,$; le pr\'e-sous-groupe $\;\hat e\;$ de $\;W\;$ est
un sous-groupe (resp. un co-cous-groupe) si et seulement si $\;f\;$ (resp.
$\;g\,$) est un sous-groupe de $\;V\,$.}

\pf  a) Comme $\;e\;$ est proportionnel \`a $\;\lambda _fg\;$ et $\;\rho
_gf\,,\;\lambda _f\in A_{f,f}\;$ et $\;\rho _g\in B_{g,g}\,$, la premi\`ere
assertion r\'esulte du corollaire 3.12\ref. Notons les co\"id\'eaux relatifs
\`a $\;W\;$ avec un exposant $\;W\,$. On a $\;G^W_e=\{\,x\in A_{f,f}\,,\;
[x,(\id\otimes \omega_{e,e}) (W)]=0\,\}\supset U\widehat D_fU\,$, vu que
$\;(\id\otimes \omega_{e,e}) (W)=\rho_g\,$.  Comme
$\;\dim\,G^W_e=\langle g,e\rangle ^{-2}=\dim \,\widehat D_f\,$, on a
l'\'egalit\'e. De m\^eme $\;\widehat D^W_e=\widehat G_g\,$.

\noindent L'assertion b) se d\'emontre de fa\c con analogue.

\noindent c) $\;e\;$ est un sous-groupe de $\;W\;$ si et seulement si
$\;U(\omega_{e,e}\otimes \id)(W)U=(\omega_{e,e}\otimes \id)(W)\,$; par
a), cela a lieu si et seulement si $\;f\;$ est un co-sous-groupe de $\;V\,$.
Les autres assertions se montrent de fa\c con analogue.

\cqfd

\rem On d\'eduit de la prop. \sn13 que l'alg\`ebre $\;A_{e,e}^W\;$ associ\'ee
\`a l'unitaire multiplicatif $\;W\;$ et \`a son pr\'e-sous-groupe $\;e\;$  est
$\;U\widehat SU\,$. De m\^eme, l'alg\`ebre $\;B_{\hat e,\hat e}^W\;$ est
$\;S\,$.

 La construction ci dessus associe \`a un triplet $\;(V,f,g)\;$
consistant en un unitaire multiplicatif $\;V\;$ et deux pr\'e-sous-groupes
$\;f\;$ et $\;g\;$ le plus \'eloign\'es possible, un unitaire multiplicatif
$\;W\,$. Si on applique \`a
nouveau cette construction au triplet $\;(W,e,\hat e)\,$, on obtient donc
$\;\widehat V\,$.

Si dans la construction de l'unitaire multiplicatif $\;W\;$ on
\'echange les r\^oles de $\;f\;$ et $\;g\,$, on obtient un unitaire
multiplicatif $\;W'\,$, de vecteur fixe $\;f\,,$ cofixe $\;g\,$, les alg\`ebres
$\;S\;$ et $\;\widehat S\;$ associ\'ees \'etant $\;A_{g,g}\;$ et $\;B_{f,f}\,$.
Or $\;A_{g,g}=\widehat JB_{g,g}\widehat J\;$ et $\;B_{f,f}=\widehat
JA_{f,f}\widehat J\,$; donc l'unitaire multiplicatif associ\'e est
$\;(\widehat J\otimes \widehat J)\Sigma W^*\Sigma (\widehat J\otimes
\widehat J)=(\widehat JJ_B\otimes\widehat JJ_B)\widehat W
(J_B\widehat J\otimes  J_B\widehat J)\,$, qui est \'equivalent \`a
$\;\widehat W\,$.

Si on part de $\;(W,\hat e,e)\,$, on obtient donc l'unitaire
$\;V'=(J_B\widehat J\otimes J_B\widehat J)V(\widehat
JJ_B\otimes \widehat JJ_B)\;$ \'equivalent
\`a $\;V\,$.

\medskip Plus g\'en\'eralement, soient $\;f,\hat f,g,\hat g\;$ des
pr\'e-sous-groupes
tels que $\;\hat f\prec f\,,\;\hat g\prec g\;$ et $\;\langle f,\hat
g\rangle=\langle
g,\hat f\rangle=n^{-1/2}\,$. Notons $\;W\;$ l'unitaire multiplicatif
associ\'e \`a
$\;(V,\hat f,g)\,$. On v\'erifie alors que:

\noindent -- $\;f\;$ et $\;\hat g\;$ sont des pr\'e-sous-groupes de $\;W\,$.

\noindent -- Les unitaires multiplicatifs construits \`a partir de
$\;(V,f,\hat g)\;$ et
$\;(W,f,\hat g)\;$ co\"incident.

On en d\'eduit alors que $\;\hat f \;$ et $\;g\;$ sont des pr\'e-sous-groupes de
l'unitaire multiplicatif $\;W'\;$ construit \`a partir de $\;(V,f,\hat
g)\;$ et que
l'unitaire multiplicatif construit \`a partir de $\;(W',\hat f,g)\;$ est
conjugu\'e \`a $\;W\,$.
\sec{Liens avec les sous-facteurs.}

\subsec{Inclusions de profondeur 2.}

Soit $\;N\subset M\;$ une inclusion irr\'eductible de profondeur
$\;2\;$ d'indice fini de facteurs de von Neumann de type $\;II_1\,$.
\hfill\break
 --- Notons $\;N\subset M\subset M_1\subset M_2\ldots\;$ la tour de Jones
de l'inclusion $\;N\subset M\,$. \hfill\break
 --- Posons $\;\widehat S=M_1\cap N'\;$ et $\;S=M_2\cap M'\;$ contenues
dans $\;M_2\cap N'\;$ qui est isomorphe \`a un $\;\L(\H)\;$ (par
d\'efinition de la profondeur 2). \hfill\break
 --- Notons $\;q\in M_1\cap N'\;$ le projecteur de Jones de l'inclusion
$\;N\subset M\;$ et $\;p\in M_2\cap M'\;$ le projecteur de Jones de
l'inclusion $\;M\subset M_1\,$. Ce sont des projecteurs minimaux de
$\;M_2\cap N'=\L(\H)\,$. Choisissons $\;e,\hat e\in \H\;$ de norme
$\;1\;$ tels que $\;p\hat e=\hat e\,,\;qe=e\;$ et $\;\langle e,\hat
e\rangle \in \R_+\,$.

Notons $\;\tau \;$ la trace normalis\'ee de $\;M_2\,$. Remarquons que,
pour $\;x\in S\,$, on a $\;E_{M_1}(x)\in M_1\cap M'=\C1\,$; donc
$\;E_{M_1}(x)=\tau (x)\,$. Donc, pour $\;x\in S\;$ et $\;y\in M_1\;$ on a
$\;\tau (xy)=\tau (x)\tau (y)\,$.

Notons alors $\;V\in \L(\H\otimes \H)\;$ l'unitaire tel que, pour
$\;a,x\in S\;$ et $\;b,y\in \widehat S\;$ on ait $\;\langle xe\otimes y\hat
e,V(ae\otimes b\hat e)\rangle =\tau (x^*y^*ab)\,$.

Notons aussi $\;W\in \L(\H\otimes L^2(M_1))\;$  tel que, pour
$\;a,x\in S\;$ et $\;b,y\in M_1\;$ on ait $\;\langle xe\otimes y\xi
_\tau ,W(ae\otimes b\xi _\tau )\rangle =\tau (x^*y^*ab)\,$. C'est l'image
de $\;V\;$ par la suite des inclusions $\;\L(\H)\otimes \L(\H)\subset
\L(\H)\otimes M_2\subset \L(\H\otimes L^2(M_1))\,$.

Il est clair que $\;W\;$ commute \`a $\;1\otimes b\;$ pour tout $\;b\in
M\,$, d'o\`u l'on d\'eduit que $\;W\in \L(\H)\otimes S\,$. Par le
th\'eor\`eme 2.\ref8 (en utilisant le lemme 2.2.c\ref), on en d\'eduit que
$\;S\;$ et $\;\widehat S\;$ sont des alg\`ebres de Hopf en dualit\'e. De
plus, pour $\;x,y\in M\;$ on a $\;W^*(1\otimes xqy)W=(1\otimes
x)W^*(1\otimes q)W(1\otimes y)\in \widehat S\otimes M_1\,$. On en
d\'eduit que l'application $\;z\mt W^*(1\otimes z)W\;$ est une coaction
(\`a gauche) de $\;\widehat S \;$ dans $\;M_1\;$ dont l'alg\`ebre des points
fixes est $\;M\,$.

On retrouve ainsi un r\'esultat d'Ocneanu; voir [\ref\Da, \Sz] pour
d'autres d\'emonstrations (dans le cadre des inclusions d'indice fini) et
[\ref\Longo, \EN] (dans un cadre plus g\'en\'eral).

\subsec{Pr\'e-sous-groupes et facteurs interm\'ediaires.}

Conservons les notations ci-dessus. Soit $\;f\;$ un pr\'e-sous-groupe de
$\;V\,$. Posons $\;M_f=\{\,x\in M_1\,,\;[x,L _f]=0\,\}\,$. C'est un
sous-facteur de $\;M_1\;$ contenant $\;M\;$ et $\;\widehat D_f\subset
\widehat S=M_1\cap N'\,$. On en d\'eduit que $\;[M_f:M]\ge \dim
\widehat D_f\,$. Par ailleurs, comme $\;L_f\ge p\,$, on a $\;L_f\xi
_\tau =\xi _\tau \,$; donc, pour $\;x\in M_f\,$, on a $\;L_fx\xi _\tau
=xL_f\xi _\tau =x\xi _\tau \,$; on en d\'eduit que $\;[M_f:M]\le \tau
(L_f)/\tau (p)=\dim \widehat D_f\,$; on en d\'eduit que $\;L_f\;$est le
projecteur de Jones de $\;M_f\subset M_1\,$.

R\'eciproquement, soit $\;M\subset P\subset M_1\;$ un facteur
interm\'ediaire. Posons $\;Q=J_MP'J_M\,\subset\L (L^2(M))\,$; on a
$\;N\subset Q\subset
M\,$; notons $\;p_1\in S\;$ le projecteur de Jones de l'inclusion
$\;P\subset M_1\;$ et $\;q_1\in \widehat S\;$ le projecteur de Jones de
l'inclusion $\;Q\subset M\,$. Comme $\;[P:M]=[M:Q]\,$, on a
$\;[P:M][Q:N]=[M:N]\,$, donc $\;\tau (p_1q_1)=\tau (p_1)\tau
(q_1)=\tau (p)[P:M]\tau (q)[Q:N]=\dim(\H)^{-1}\,$. Par ailleurs, $\;q_1\in
P\;$ (car $\;P\;$ est la construction de base de $\;Q\subset M\;$ et
$\;p_1\in P'\,$; donc $\;p_1q_1=q_1p_1\,$. Par
le corollaire 3.12\ref, il existe un pr\'e-sous-groupe $\;f\;$ tel que
$\;p_1=L_f\;$ et $\;q_1=\rho _f\,$.

En d'autres termes, les $\;L_f\;$ o\`u $\;f\;$ est un pr\'e-sous-groupe,
sont les projecteurs de Jones des facteurs interm\'ediaires $\;M\subset
P\subset M_1\,$. L'application $\;P\mt J_MP'J_M\;$ est une
correspondance bijective entre facteurs interm\'ediaires $\;M\subset
P\subset M_1\;$ et facteurs interm\'ediaires $\;N\subset
Q\subset M\,$. On en d\'eduit que les $\;\rho _f\;$ o\`u $\;f\;$ est un
pr\'e-sous-groupe, sont les projecteurs de Jones des facteurs
interm\'ediaires $\;N\subset Q\subset M\,$.

On en d\'eduit une bijection (d\'ecroissante) $\;f\mt N_f\;$ entre
pr\'e-sous-groupes et facteurs interm\'ediaires  $\;N\subset P\subset M\,$.
On retrouve ainsi, \`a l'aide de la prop. 4.3\ref, une bijection croissante (un
isomorphisme d'ensembles ordonn\'es) entre  facteurs interm\'ediaires et
sous-alg\`ebres co\"id\'eales \`a gauche de l'alg\`ebre de Hopf $\;S\;$
[\ref\ILP] (voir aussi [\Enockb] o\`u ce r\'esultat est g\'en\'eralis\'e au cas
non compact).

\medskip Plusieurs de nos constructions s'interpr\^etent alors:

\noindent -- Le th\'eor\`eme de finitude (corollaire 3.6\ref), est donc une
cons\'equence du fait que dans une inclusion d'indice fini il y a un nombre
fini de facteurs interm\'ediaires ([\ref\Wa]).

\noindent -- Les sous-groupes correspondent aux facteurs
interm\'ediaires $\;N_f\;$ tels que l'inclusion $\;N_f\subset M\;$ soit de
profondeur $\;2\,$; les co-sous-groupes correspondent aux facteurs
interm\'ediaires $\;N_f\;$ tels que l'inclusion $\;N\subset N_f\;$ soit de
profondeur $\;2\,$.

\noindent -- Les sous-quotients correspondent aux couples de facteurs
interm\'ediaires $\;N_{\hat f}\,,\; N_f\;$ tels que $\;N_f\subset N_{\hat
f}\;$ et que l'inclusion $\;N_f\subset N_{\hat f}\;$ soit de profondeur
$\;2\,$.

\noindent -- Des pr\'e-sous-groupes $\;f\;$ et $\;g\;$ sont le plus
\'eloign\'es possible si et seulement si on a un carr\'e commutatif et
cocommutatif
$$\matrix{N&\subset &N_f\,\cr \cap&&\cap\cr
N_g& \subset &M\,.}$$Notons alors $\;M_f\;$ la construction de base de
Jones de l'inclusion $\;N_f\subset M\,$.  La construction du paragraphe 6,
dit alors que l'inclusion $\;N_g\subset M_f\;$ est irr\'eductible de
profondeur $\;2\,$.

On peut en fait donner une d\'emonstration directe de ce r\'esultat:

\exer On a $\;N_g'\cap M_1=\widehat G_g\subset \widehat S=N'\cap
M_1\,$. On a $\;N'\cap M_f=\widehat D_f\,$. Donc $\;N_g'\cap
M_f=(N_g'\cap M_1)\cap (N'\cap M_f)=\widehat G_g\cap \widehat
D_f=\C\,$.

\exer La construction de base appliqu\'ee \`a l'inclusion $\;N_g\subset
M_f\;$ donne un facteur $\;M^1_g\,$. On a $\;M'\cap M^1_g=D_g\subset
S=M'\cap M_2\,$. On en d\'eduit que $\;N_g'\cap M^1_g\;$ contient
l'alg\`ebre $\;B_{g,g}\,$. Sa dimension est alors \'egale \`a l'indice de
$\;N_g\subset M_f\,$, d'o\`u le r\'esultat.

\bigskip\bigskip\centerline{\twelvebf BIBLIOGRAPHIE.}

[\BSb] S. Baaj et G. Skandalis: Unitaires multiplicatifs et dualit\'e pour
les produits crois\'es de $\,C^*$-alg\`ebres. Annales Scient. E. Norm.
Sup. 4e s\'erie {\bf 26} (1993) 425-488.

[\Da] M.-C. David: Paragroupe d'Adrian Ocneanu et alg\`ebre de Kac. Pac.
J. of Math. {\bf 172} No 2 (1996), 331-363.

[\Dab] M.-C. David: Couple assorti de syst\`emes de Kac et inclusions de
facteurs de type $II_1$. {\it Pr\'epubli\-cation.}

[\Enock] M. Enock: Inclusions irr\'eductibles de facteurs et unitaires
multiplicatifs II. {\it Pr\'epublica\-tion.}

[\Enockb] M. Enock: Sous-facteurs interm\'ediaires et groupes
quantiques mesur\'es. {\it Pr\'epublica\-tion.}

[\EN] M. Enock and R. Nest: Irreducible inclusions of factors,
multiplicative unitaries and Kac algebras. J.F.A. {\bf 137}, No 2 (1996),
466-543.

[\ES] M. Enock et J.-M. Schwartz: {\it Kac Algebras and Duality of
Locally Compact Groups.}\/ Springer (1992).

[\GHJ] F.M. Goodman, P. de la Harpe, V.F.R. Jones: {\it Coxeter graphs and
towers of algebras.}\/ M.S.R.I. publ. {\bf 14}.

[\HS] J.H. Hong and W. Szyma\'nski: Composition of subfactors and
twisted  bicrossed products. J. Operator Theory 37 (1997), no. 2,
281-302.

[\IK] M. Izumi and H. Kosaki: Finite-dimensional Kac algebras arising
from  certain group actions on a factor. Internat. Math. Res. Notices
1996, no. 8,  357-370.

[\ILP] M. Izumi, R. Longo and S. Popa: A Galois Correspondance for
Compact Groups of Automorphisms of von Neumann Algebras with a
Generalization to Kac Algebras. Preprint Feb. 1996.

[\Kac] G.I. Kac: Extensions of groups to ring groups. Math U.S.S.R. Sbornik
{\bf 5} (1968), 451-474.

[\KP] G.I. Kac and V.G. Paljutkin: Finite group rings. Trans. Moskow Math.
Soc. (1966), 251-294.

[\Longo] R. Longo: A duality for Hopf algebras and subfactors I. Comm.
Math. Phys. {\bf 159} (1994), 133-155.

[\Maj] S. Majid: Hopf-von Neumann algebra bicrossproducts, Kac algebra
bicrossproducts, and the classical Yang-Baxter equations. J.F.A. {\bf
95} No 2 (1991), 291-319.

[\Ma] A. Masuoka: Freeness of Hopf Algebras over coideal subalgebras.
Comm. in Alg. {\bf 20} (5) (1992), 1353-1373.

[\Mon] S. Montgomery: {\it Hopf algebras and their action on rings.}\/
C.B.M.S. lecture notes, {\bf 82}, A.M.S. (1993).

[\NZ] W.D. Nichols and M.B. Zoeller: Hopf algebra freeness theorem. Amer.
J. of Math. {\bf 111} (1989), 381-385.

[\Po] S. Popa: Orthogonal pairs of $*$-subalgebras in finite von Neumann
algebras. J. Operator Theory {\bf 9} (1983), 253--268.

[\Sz] W. Szymanski: Finite index subfactors and Hopf algebras crossed
products. Proc. A.M.S. {\bf 120} (1994), 519-528.

[\Takeu] M. Takeuchi: Matched pairs of groups and bismashed product of
Hopf algebras. Comm. Alg. {\bf 9} (1981), 841-882.

[\Wa] Y. Watatani: Lattice structure of intermediate subfactors. Math.
Phys. Stud., {\bf 16}, Quantum and non-commutative
analysis (Kyoto, 1992), (1993), 331-333.

\bye